\documentclass[11pt]{article}

\usepackage{amsfonts, amssymb}
\usepackage{amsmath}
\usepackage{amsthm}
\usepackage{latexsym}
\usepackage{verbatim}
\usepackage{psfig}
\usepackage{url}
\usepackage{graphicx,color}
\definecolor{lightgray}{gray}{.90}

\def\t{\mbox{\textbf{\textsf{t}}}}

\newcommand{\NN}{\mathbb{N}}

\newcommand{\ZZ}{\mathbb{Z}}
\newcommand{\RR}{\mathbb{R}}
\newcommand{\TT}{\mathrm{T}}

\def\Alph{\textrm{Alph}}
\def\cE{\mathcal{E}}
\def\cS{\mathcal{S}}
\def\cP{\mathcal{P}}

\def\br{\mathbf{r}}
\def\bw{\mathbf{w}}

\advance \topmargin by -\headheight
\advance \topmargin by -\headsep
\newenvironment{debut}[1]%
{\begin{quote} \bf #1 \it}%
{\end{quote}}

{\begin{debut} {R\'esum\'e~: }}
{\end{debut}}

{\begin{debut} {Abstract:}}
{\end{debut}}

\theoremstyle{plain}
\newtheorem{theo}{Theorem}[section]

\newtheorem{cor}[theo]{Corollary}
\newtheorem{prop}[theo]{Proposition}
\newtheorem{defi}[theo]{Definition}

\theoremstyle{definition}
\newtheorem{remq}[theo]{Remark}
\newenvironment{rem}{\begin{remq}\rm }{\end{remq}}
\newtheorem{exa}[theo]{Example}
\newtheorem*{notation}{Notation}

\theoremstyle{remark}
\newtheorem*{note}{Note}

\def\t1{\mathbf{t}^{(1)}}
\def\bt{\mathbf{t}}
\def\bm{\mathbf{m}}
\def\m1{\mathbf{m}^{(1)}}

\def\bs{\mathbf{s}}
\def\cA{\mathcal{A}}
\def\bbf{\mathbf{f}}

\def\bb{\mathbf{b}}
\def\bu{\mathbf{u}}
\def\bv{\mathbf{v}}
\def\bx{\mathbf{x}}
\def\by{\mathbf{y}}

\def\br{\mathbf{r}}
\def\bl{\mathbf{l}}
\def\rev{\widetilde}
\def\cAstar{\mathcal{A}^*}

\def\cAw{\mathcal{A}^\omega}
\def\cAinf{\mathcal{A}^\infty}
\def\cB{\mathcal{B}}

\def\empt{\varepsilon}
\def\Ult{\mathrm{Ult}}
\def\Alphit{\mbox{{\em Alph}}}

\numberwithin{equation}{section}

\begin{document}

\title{\vspace{-0cm} EPISTURMIAN WORDS: A SURVEY\footnote{This paper grew out of an invited lecture given by the second author at the {Sixth International Conference on Words}, Marseille, France, September 17--21, 2007.}}
\author{Amy Glen\footnotemark[2] \and Jacques Justin\footnotemark[3]}  
\date{Submitted: December 11, 2007; Revised: September 16, 2008}          
\maketitle  

\footnotetext[2]{Corresponding author: LaCIM, Universit\'e du Qu\'ebec \`a Montr\'eal, C.P. 8888, succursale Centre-ville, Montr\'eal, Qu\'ebec, H3C 3P8, CANADA $\backslash$ The Mathematics Institute, Reykjavik University, Kringlan 1, IS-103 Reykjavik, ICELAND (\texttt{amy.glen@gmail.com}).}
\footnotetext[3]{LIAFA, Universit\'e Paris  Diderot - Paris 7,  
Case 7014, 75205 Paris Cedex~13, FRANCE    
(\texttt{jacjustin@free.fr}).
}


\hrule
\begin{abstract} 
In this paper, we survey the rich theory of infinite {\em episturmian words} which generalize to any finite alphabet, in a rather resembling way, the well-known family of {\em Sturmian words} on two letters. After recalling definitions and basic properties, we consider {\em episturmian morphisms} that allow for a deeper study of these words. 
Some properties of factors are described, including factor complexity, palindromes, fractional powers, frequencies, and return words. We also consider lexicographical properties of episturmian words, as well as their connection to the balance property, and related notions such as finite episturmian words, Arnoux-Rauzy sequences,  and ``episkew words'' that generalize the skew words of Morse and Hedlund. \medskip

\medskip

\end{abstract}
\hrule

\section{Introduction} 

\subsection{From Sturmian to episturmian}

Most renowned amongst the branches of combinatorics on words is the theory of infinite binary sequences called \emph{Sturmian words}, which are fascinating in many respects, having been studied from combinatorial, algebraic, and geometric points of view. Their beautiful properties are related to many fields such as Number Theory, Geometry, Symbolic Dynamical Systems, Theoretical Physics, and Theoretical Computer Science (see \cite{jAjS03auto, mL02alge, nP02subs} for recent surveys).



Since the seminal works of Morse and Hedlund~\cite{gHmM40symb}, Sturmian words have been shown to admit numerous equivalent definitions and characterizations. 
For instance, it is well known that an infinite word $\bw$ over $\{a,b\}$ is Sturmian if and only if $\bw$ is aperiodic and {\em balanced}: for any two factors $u$, $v$ of $\bw$ of the same length, the number of $a$'s  in each of $u$ and $v$ differs by at most $1$.  Sturmian words are also characterized by their {\em factor complexity} function (which counts the number of distinct factors of each length): they have exactly $n + 1$ distinct factors of length $n$ for each $n$. In this sense, Sturmian words are precisely the aperiodic infinite words of minimal factor complexity since, as is well known, an infinite word is ultimately periodic if and only if it has less than $n+1$ factors of length $n$ for some $n$ (see \cite{eCgH73sequ}). Many interesting properties of Sturmian words can be attributed to their low complexity, which induces certain regularities in such words without, however, making them periodic. 
Sturmian words can also be geometrically realized as {\em cutting sequences} by considering the sequence of `cuts' in an integer grid made by a line of irrational slope (see for instance \cite{dCwMaPpS93subs, yB95comp}).  They also provide a symbolic coding of the orbit of a point on a circle with respect to a rotation by an irrational number (see \cite{gHmM40symb, pAvB98three}). 

All of the above characteristic properties of Sturmian words lead to natural generalizations on arbitrary finite alphabets. In one direction, the balance property naturally extends to an alphabet with more than two letters (e.g., see \cite{pH00suit,rT00frae, lv03bala}) as does the following {\em generalized balance property} that also characterizes Sturmian words (see \cite{iFlV02gene, bA03bala}): the difference between the number of occurrences of a word $u$ in any pair of factors of the same length is at most $1$. In another direction, we could consider relaxing the minimality condition for the factor complexity $p(n)$. For example, {\em quasi-Sturmian words} are infinite words for which there exist two positive integers $N$ and $c$ such that $n+1 \leq p(n) \leq n + c$ for all $n\geq N$. This generalization was introduced in \cite{jAjDmQlZ01tran} when studying the transcendence of certain continued fraction expansions. See also  \cite{jC98sequ, eC74sequ, aHrT00char, rRlZ00agen} for similar extensions of Sturmian words with respect to factor complexity. From the geometric point of view, cutting sequences naturally generalize to trajectories in the hypercube billiard (e.g., see \cite{jBcR05pali}), and codings of rotational orbits carry over to codings of interval exchange transformations (e.g., see \cite{jBlV02codi}). 

Two other very interesting natural generalizations of Sturmian words are {\em Arnoux-Rauzy sequences} \cite{pAgR91repr, gR82suit} and {\em episturmian words} \cite{xDjJgP01epis,jJgP02epis}, which we will now define. 

From the factor complexity of Sturmian words, it immediately follows that any Sturmian word is over a $2$-letter alphabet and has exactly one {\em left special} factor of each length. A factor $u$ of a finite or infinite word $w$ is said to be {\em left special} (resp.~{\em right special}) in $w$ if there exists at least two distinct letters $a$, $b$ such that $au$ and $bu$ (resp.~$ua$, $ub$) are factors of $w$. Extending the left special property of Sturmian words, a recurrent infinite word $\bw$ over a finite alphabet $\cA$ is said to be an {\em Arnoux-Rauzy sequence} (or a {\em strict episturmian word}) if it has exactly one left special factor and one right special factor of each length, and for every left (resp.~right) special factor $u$ of $\bw$, $xu$ (resp.~$ux$) is a factor of $\bw$ for all letters $x \in \cA$. A noteable property that is shared by Sturmian words and Arnouxy-Rauzy sequences is their closure under reversal, i.e., if $u$ is a factor of such a word, then its reversal is also a factor. This nice property inspired Droubay, Justin, and Pirillo's generalization of Sturmian words in~\cite{xDjJgP01epis}: an infinite word is {\em episturmian} if it is closed under reversal and has at most one left special factor of each length. Sturmian, Arnoux-Rauzy, and episturmian words all have {\em standard} (or {\em characteristic}) elements, which are those having all of their left special factors as prefixes. Within these families of words, standard words are good representatives in the sense that an infinite word belongs to one such family if and only if it has the same set of factors as some standard word in that family.

From the definitions, it is clear that the family of Arnoux-Rauzy sequences is a particular subclass of the family of episturmian words. More precisely, episturmian words are composed of  the Arnoux-Rauzy sequences, images of the Arnoux-Rauzy sequences by {\em episturmian morphisms}, and certain periodic infinite words (see Section~\ref{S:A-R}). In the $2$-letter case, Arnoux-Rauzy sequences are exactly the Sturmian words whereas episturmian words include all recurrent balanced words, i.e., periodic balanced words and Sturmian words.

The study of episturmian words and Arnoux-Rauzy sequences has enjoyed a great deal of popularity in recent times, owing mostly to the many properties that they share with Sturmian words. 
In this paper we survey the purely combinatorial work on episturmian words,  beginning with their definition and basic properties in Section \ref{S:definition}. Then, in Section~\ref{S:epi-morphisms}, we recall episturmian morphisms which allow for a deeper study of episturmian words. In particular, any episturmian word is the image of another episturmian word by some so-called {\em pure episturmian morphism}. Even more, any episturmian word can be infinitely decomposed over the set of pure episturmian morphisms.  This last property allows  an episturmian word to be defined by one of its morphic decompositions or, equivalently, by a certain {\em directive word}, which is an infinite sequence of rules  for decomposing the given episturmian word by morphisms. In Section~\ref{S:shifts&spins} we consider notions such as {\em shifts}, {\em spins}, and {\em block-equivalence} in connection with directive words, which allow us to study when two different {\em spinned infinite words} direct the same episturmian word. We also consider periodic and purely morphic episturmian words. In Section~\ref{S:A-R}, our discussion briefly turns to Arnoux-Rauzy sequences and {\em finite episturmian words}. Following this, we study in Section~\ref{S:factors} some properties of factors of episturmian words (and Arnoux-Rauzy sequences), including factor complexity, palindromes, fractional powers, frequencies, and return words.  Lastly, we consider more recent work involving {\em lexicographic order} and the balance property (including {\em Fraenkel's conjecture}). 

\subsection{Notation \& terminology}

We assume the reader is familiar with combinatorics on words and morphisms (e.g., see \cite{mL83comb,mL02alge}). In this section, we recall some basic definitions and properties relating to episturmian words which are  needed throughout the paper.  For the most part, we follow the notation and terminology of \cite{xDjJgP01epis, jJgP02epis, jJgP04epis, aGjJgP06char}.

Let $\cA$ denote a finite {\em alphabet}, i.e., a non-empty finite set of symbols called {\em letters}. A finite \emph{word} over $\cA$ is a finite sequence of letters from $\cA$. The {\em empty word} $\empt$ is the empty sequence. Under the operation of concatenation, the set $\cA^*$ of all finite words over $\cA$ is a {\em free monoid} with identity element $\empt$ and set of generators $\cA$. The set of {\em non-empty} words over $\cA$ is the {\em free semigroup} $\cA^+ := \cAstar \setminus \{\empt\}$.  

A {\em right-infinite} (resp.~{\em left-infinite}, {\em bi-infinite}) word over $\cA$ is a sequence indexed by $\NN^+$ (resp.~$\ZZ\setminus\NN^+$, $\ZZ$) with values in $\cA$. For instance, a left-infinite word is represented by $\bu = \cdots b_{-2}b_{-1}b_0$ and a right-infinite word by $\bv = b_1b_2b_3\cdots$ where $b_i \in \cA$. The concatenation of $\bu$ and $\bv$ gives the bi-infinite word $\bu.\bv = \cdots b_{-2}b_{-1}b_{0}.b_1b_2b_3\cdots$ with a dot written between $b_0$ and $b_1$ to avoid ambiguity. For easier reading, infinite words are hereafter typically typed in boldface to distinguish them from finite words.


The {\em shift map} $\TT$ is defined for bi-infinite words $\mathbf{b} = (b_i)_{i\in\ZZ}$ by $\TT(\mathbf{b}) = (b_{i+1})_{i\in\ZZ}$ and its $k$-th iteration is denoted by $\TT^k$. This extends to right-infinite words for $k\geq 0$ and left-infinite words for $k\leq 0$. For finite words $w \in \cA^*$, the shift map $\TT$ acts circularly, i.e., if $w = xv$ where $x \in \cA$, then $\TT(w) = vx$.

The set of all right-infinite words over $\cA$ is denoted by $\cAw$, and we define $\cAinf := \cAstar \cup \cAw$.  An {\em ultimately periodic} right-infinite word can be written as $uv^\omega = uvvv\cdots$, for some $u$, $v \in \cAstar$, $v\ne \empt$. If $u = \empt$, then such a word is {\em periodic}. A right-infinite word that is not ultimately periodic is said to be {\em aperiodic}.

Given a finite word $w = x_{1}x_{2}\cdots x_{m} \in \cA^*$ with each $x_{i} \in \cA$, the \emph{length} of $w$, denoted by $|w|$, is equal to $m$. By convention, the empty word $\empt$ is the unique word of length $0$.   The number of occurrences of a letter $a$ in $w$ is denoted by $|w|_a$. If $|w|_a = 0$, then $w$ is said to be {\em $a$-free}.  The \emph{reversal} $\rev{w}$ of $w$ is its mirror image: $\rev{w} = x_{m}x_{m-1}\cdots x_{1}$, and if $w = \rev{w}$, then $w$ is called a \emph{palindrome}. The reversal operator naturally extends to bi-infinite words; that is, the reversal of the bi-infinite word $\mathbf{b} = \bl.\br$, with $\bl$ left-infinite and $\br$ right-infinite, is given by $\rev{\mathbf{b}} = \rev\br.\rev\bl$.

A finite word $w$ is a \emph{factor} of a finite or infinite word $z$ if $z = uwv$ for some words $u$, $v$ (which are finite or infinite depending on $z$). In the special case $u = \empt$ (resp.~$v =\empt$),  we call $w$ a \emph{prefix} (resp.~\emph{suffix}) of $z$.  We use the notation $p^{-1}w$ (resp.~$ws^{-1}$) to indicate the removal of a prefix $p$ (resp.~suffix $s$) of a finite word $w$. Note that a prefix or suffix $u$ of a finite word $w$ is said to be {\em proper} if $u \ne w$. A factor $u$ of a finite or infinite word $w$ is \emph{right} (resp.~\emph{left}) \emph{special} if $ua$, $ub$ (resp.~$au$, $bu$) are factors of $w$ for some letters $a$, $b \in \cA$, $a \ne b$. 

For any finite or infinite word $w$, $F(w)$ denotes the set of all its factors. Moreover, the \emph{alphabet} of $w$ is Alph$(w):= F(w) \cap \cA$ and, if $w$ is infinite, we denote by Ult$(w)$ the set of all letters occurring infinitely often in $w$.  Any two infinite words $\bx$, $\by$ are said to be \emph{factor-equivalent} if $F(\bx) = F(\by)$, i.e., if $\bx$ and $\by$ have the same set of factors. 

A factor of an infinite word $\bx$ is \emph{recurrent} in $\bx$ if it occurs infinitely often in $\bx$, and $\bx$ itself is said to be \emph{recurrent} if all of its factors are recurrent in it. For a bi-infinite word to be recurrent, any factor must occur infinitely often to the left and to the right. An infinite word is said to be \emph{uniformly recurrent} if any factor occurs infinitely many times in it with bounded gaps~\cite{eCgH73sequ}.

A \textit{morphism} $\varphi$ on $\cA$ is a map from $\cA^*$ to $\cA^*$ such that $\varphi(uv) = \varphi(u) \varphi(v)$ for any words $u$, $v$ over $\cA$.  A morphism on $\cA$ is entirely defined by the images of letters in $\cA$. All morphisms considered in this paper will be {\em non-erasing}: the image of any non-empty word is never empty. Hence the action of a morphism $\varphi$ on $\cA^*$ can be naturally extended to infinite words; that is, if $\bx = x_1x_2x_3 \cdots \in \cAw$, then $f(\bx) = f(x_1)f(x_2)f(x_3)\cdots$. An infinite word $\bx$ can therefore be a {\em fixed point} of a morphism $\varphi$, i.e., $\varphi(\bx) = \bx$.  If $\varphi$ is a (non-erasing) morphism such that 
$\varphi(a) = aw$ for some letter $a \in \cA$ and $w \in \cA^{+}$,
then $\varphi^{n}(a)$ is a proper prefix of the word $\varphi^{n+1}(a)$ for each $n \in \NN$, and the limit of the sequence $(\varphi^n(a))_{n\geq 0}$ is the unique infinite word: 
\[
  \bw = \underset{n \rightarrow\infty}{\lim}\varphi^n(a) = \varphi^{\omega}(a) ~(= aw\varphi(w)\varphi^2(w)\varphi^3(w)\cdots).
\]
Clearly, $\bw$ is a fixed point of $\varphi$ and we say that $\bw$ is \emph{generated} by $\varphi$. Furthermore, an infinite word generated by a morphism is said to be {\em purely morphic}. 

In what follows, we will denote the composition of morphisms by juxtaposition as for concatenation of words.

\section{Definitions \&  basic properties} \label{S:definition}

In the initiating paper \cite{xDjJgP01epis}, episturmian words were defined as an extension of {\em standard episturmian words}, which were first introduced as a generalization of {\em standard} (or {\em characteristic}) Sturmian words using iterated palindromic closure (a construction due to de~Luca \cite{aD97stur}). Here we choose instead to begin with the following definition for deriving the main basic properties of episturmian words.

\begin{defi} \label{D:main} {\em \cite{xDjJgP01epis}}
An infinite word $\bt \in \cAw$ is \emph{episturmian} if $F(\bt)$ is closed under reversal and $\bt$ has at most one left special factor (or equivalently, right special factor) of each length. Moreover, an episturmian word is \emph{standard} if all of its left special factors are prefixes of it.
\end{defi}
\begin{note} We can equivalently consider left or right special factors in the first part of the above definition since, by closure under reversal, a factor is left (resp.~right) special if and only if its reversal is right (resp.~left) special.
\end{note}

\begin{rem} When $|\cA| = 2$, Definition~\ref{D:main} gives the (aperiodic) Sturmian words, as well as the periodic balanced infinite words (also known as the {\em periodic Sturmian words}). See for instance~\cite{aGjJgP06char} or Section~\ref{SS:q-balance}.
\end{rem}

The following theorem collects together some useful characteristic properties of standard episturmian words. Before stating it, let us first recall the some definitions. 

Given two palindromes $p$, $q$, we say that $q$ is a {\em central factor} of $p$ if $p = w q \rev{w}$ for some $w \in \cA^*$. The \emph{palindromic right-closure} $w^{(+)}$ of a finite word $w$ is the (unique) shortest palindrome having $w$ as a prefix (see \cite{aD97stur}). That is, $w^{(+)} = wv^{-1}\rev{w}$ where $v$ is the longest palindromic suffix of $w$. For example, $(race)^{(+)} = race\thinspace car$. The {\em iterated palindromic closure} function \cite{jJ05epis}, denoted by $Pal$, is defined recursively as follows. Set $Pal(\empt) = \empt$ and, for any word $w$ and letter $x$, define $Pal(wx)~=~(Pal(w)x)^{(+)}$. For instance, $Pal(abc)=(Pal(ab)c)^{(+)} = (abac)^{(+)} = abacaba$. (See Sections~\ref{SS:pure} and~\ref{SS:palindromic-closure} for further insight about palindromic closure.)

\begin{theo} \label{T:main}
For an infinite word $\bs \in \cAw$, the following properties are equivalent.
\begin{enumerate}
	\item[i)] $\bs$  is standard episturmian.
	\item[ii)] Any first occurrence of a palindrome in $\bs$ is a central factor of some palindromic prefix of $\bs$  (property Pi).
	\item[iii)] If $w$ is a prefix of $\bs$, then $w^{(+)}$ is also a prefix of $\bs$ (property Al). 
	\item[iv)] There exists an infinite word $\Delta = x_1x_2 \cdots$ $(x_i \in \cA)$, called the {\em directive word} of $\bs$, such that $\bs = \lim_{n\rightarrow\infty}Pal(x_1\cdots x_n)$. 
\end{enumerate}
\end{theo}
\begin{rem} The palindromes $Pal(x_1 \cdots x_n)$ are very often denoted by $u_{n+1}$ in the literature (and we will sometimes use the latter notation when convenient). By construction, these palindromes are exactly the palindromic prefixes of $\bs$. Moreover, $\bs$ is uniquely determined by the directive word $\Delta$. 
\end{rem}
\begin{proof}[Proof of Theorem~$\ref{T:main}$]
$i)\Rightarrow ii)$: 
Let $\bs=up\bt$, $u \in \cAstar,\,\bt \in\cAw $ showing the first occurrence of some palindrome $p$ in  $\bs$. Suppose $p$ is not the central factor of a palindromic prefix. Then we have $\bs=vxwp\tilde w y\bt'$, $x \neq y \in \cA$. By the reversal property, $y  w p \tilde w x \in F(s)$, thus $ wp \tilde w$ is left special, hence is a prefix of $\bs$. Thus $p$ has another occurrence strictly on the left of the considered one, a contradiction. \par

$i)\Rightarrow iii)$: If $iii)$ is false, let $w=ux$, with $u \in \cAstar$ and $x \in \cA$, be the shortest prefix of $\bs$ such that $w^{(+)}$ is not a prefix of $\bs$ . Thus $u^{(+)}$ is a prefix of $\bs$. If $u$ were not a palindrome then $w$ would be a prefix of $u^{(+)}$; whence $w^{(+)}=u^{(+)} $, a contradiction.
Thus $u$ is a palindrome. Now let $q$ be the longest palindromic suffix of $w$. Then 
$w^{(+)}=w_1q \tilde w_1 = w \tilde w_1$ where $w=w_1q$, and $w^{(+)}=w_1q fyg$ and $w_1qfz$ is a prefix of $\bs$  for some $y \ne z \in \cA$, $f$, $g \in \cAstar$.  
Hence $y \tilde f q \in F(\tilde w) \subset F(\bs)$ and $z \tilde f q \in F(\bs)$. Therefore $\tilde f q$ is a left special prefix of $\bs$. As $qf$ is a prefix of $\tilde w=xu$, $x^{-1} q f$ is a prefix of $u$, hence $x^{-1} q f \alpha$ is a prefix of $u$ for some letter $\alpha$. So we have $x^{-1} q f \alpha = \tilde f q$, whence $\alpha=x$ and $qfx=x \tilde f q$. This word is a palindrome and, as it is a suffix of $w$, this contradicts the minimality of $|q|$.

$iii)\Rightarrow iv)$: Trivial.

At this stage, we have proved that standard episturmian words satisfy $ii), iii),iv)$. The equivalence of these three properties is proved in \cite[Theorem 1]{xDjJgP01epis}. Finally, if $\bs$ satisfies them, then $F(\bs)$ is closed under reversal and by \cite[Proposition 5]{xDjJgP01epis} all of its left special factors  are prefixes of it, thus $\bs$ is standard episturmian. 
\end{proof}

\begin{rem}
Hereafter, we adopt ``epistandard'' as a shortcut for ``standard episturmian'', as in \cite{aGfLgR07quas, gR03conj, gR07aloc}. Also, unless stated otherwise, the notation $\Delta = x_1x_2x_3 \cdots$ ($x_i \in \cA$) will remain for the directive word of an epistandard word $\bs$.
\end{rem}
 
\begin{exa} \label{Ex:trib}
The epistandard word directed by $\Delta = (abc)^\omega$ is known as the {\em Tribonacci word} (or {\em Rauzy word} \cite{gR82suit}); it begins in the following way:
\[
  \br = \underline{a}\underline{b}a\underline{c}aba\underline{a}bacaba\underline{b}acabaabacaba\underline{c}abaabaca \cdots~,
\]  
where each palindromic prefix $Pal(x_1\cdots x_n)$ is followed by an underlined letter $x_n$. More generally, for $k\geq 2$, the {\em $k$-bonacci word} is the epistandard word over $\{a_1,\ldots, a_k\}$ directed by $(a_1a_2\cdots a_k)^\omega$ (e.g., see \cite{aG05powe}). 
\end{exa}

\begin{note}
For recent studies of the properties of Tribonacci word, see for instance  \cite{aG06onst, bTzW06some}  and the chapter by Allouche and Berth\'e in \cite{mL05appl}.
\end{note}


\subsection{Equivalence classes}

In \cite{xDjJgP01epis}, an infinite word $\bt \in \cAw$ was said to be episturmian if $F(\bt) = F(\bs)$ for some epistandard word $\bs$. This definition is equivalent to Definition \ref{D:main} by Theorem 5 in \cite{xDjJgP01epis}.  Moreover, it was proved in \cite{xDjJgP01epis} that episturmian words are uniformly recurrent, by showing that this nice property is implied by $iv)$ of Theorem \ref{T:main}. Thus,  ultimately periodic episturmian words are (purely) periodic. The aperiodic episturmian words are exactly those episturmian words with exactly one left special factor of each length.

In each {\em equivalence class} of episturmian words (i.e., same set of factors), there is one epistandard word in the aperiodic case and two in the periodic case, except if this word is $a^\omega$ with $a$ a letter. 
For example, $\bs _1= (abac)^\omega$ has directive word $\Delta _1 =abc^\omega$ and $\bs _2 =(acab)^\omega$ is directed by $\Delta _2 =acb^\omega$. Both $\bs_1$ and $\bs_2$ are standard with the same factors. Theorem~\ref{T:directSame} in Section~\ref{SS:directive} demonstrates why this is true in general (see also Remark~\ref{R:subshift}).

\subsection{Bi-infinite episturmian words} \label{SS:bi-infinite}

Definition \ref{D:main} can be extended to bi-infinite words, in which case we must assume they are recurrent. (As is well known, recurrence follows automatically from closure under reversal in the case of  right-infinite words; see for instance \cite{mBaDaGlZ08acon} for a proof of this fact.) Bi-infinite words are sometimes more natural because in particular they can be shifted in both directions, allowing for simpler formulations. More specifically, a (right-infinite) episturmian word $\bt$ can be prolonged infinitely to the left with the same set of factors, i.e., remaining in the same equivalence class. There are several or one such prolongation according to whether or not $\bt = \TT^i(\bs)$, with $\bs$ epistandard and $i \geq 0$ (see \cite{jJgP02epis, jJgP04epis}).

\begin{note}
Hereafter, `infinite word' should be taken to mean a right-infinite word, whereas left-infinite and bi-infinite words will be explicitly referred to as such. 
\end{note}

\subsection{Strict episturmian words} \label{SS:strict}

An epistandard word $\bs\in \cA^\omega$, or any factor-equivalent (episturmian) word $\bt$,    
is said to be \emph{$\cB$-strict} (or $k$-\emph{strict} if $|\cB|=k$, or {\em strict} if $\cB$ is understood) if 
Alph$(\Delta) =$ Ult$(\Delta) = \cB \subseteq \cA$.  That is, an episturmian word is strict if every letter in its alphabet occurs infinitely often in its directive word.  

The $k$-strict episturmian words are precisely the episturmian words $\bt$ having exactly one left special factor of each length and for which any left special factor $u$ in $\bt$ has $k=|\cA|$ different left extensions in $\bt$ (i.e., $xu$ is a factor of $\bt$ for all letters $x$ in the $k$-letter alphabet $\cA$). As a consequence, $k$-strict episturmian words have factor complexity $(k - 1)n + 1$ for each $n \in \NN$ (see \cite[Theorem 7]{xDjJgP01epis}); such words are exactly the $k$-letter {\em Arnoux-Rauzy  sequences}, the study of which began in \cite{pAgR91repr} (see also \cite{ jJgP02onac, rRlZ00agen} for example). In particular, the $2$-strict episturmian words correspond to the (aperiodic) Sturmian words.  Arnoux-Rauzy sequences will be discussed further in Section~\ref{S:A-R}.

\begin{rem} \label{R:periodic}
A noteworthy fact is that an episturmian word is periodic if and only if $|\Ult(\Delta)| = 1$ (see \cite[Proposition 2.9]{jJgP02epis}).  The exact form of a periodic episturmian word is given by Theorem~\ref{T:periodic} in Section~\ref{SS:rigidity}. We first need to consider {\em episturmian morphisms}.
\end{rem}


\section{Episturmian morphisms} \label{S:epi-morphisms}

From Lemma 4 in \cite{xDjJgP01epis}, if $\bs$ is epistandard with first letter $a = x_1$, then $a$ is separating for $\bs$ and its factors, i.e., any factor of $\bs$ of length $2$ contains the letter $a$. Any episturmian word $\bt$ that is factor-equivalent to $\bs$ also has separating letter $a$, and hence can be factorized with a code:
\[
\begin{cases}
\{a\} \cup a(\cA \setminus \{a\})  &\mbox{if $\bt$ begins with $a$}, \\
\{a\} \cup (\cA \setminus \{a\})a  &\mbox{otherwise}. 
\end{cases}
\] 
This leads to \emph{episturmian morphisms}, which were introduced by Justin and Pirillo \cite{jJgP02epis} in order to study deeper properties of episturmian words. As we shall see in Section~\ref{SS:relation}, episturmian morphisms are precisely the morphisms that preserve the set of aperiodic episturmian words (i.e., the morphisms that map aperiodic episturmian words onto aperiodic episturmian words). Such morphisms naturally generalize to any finite alphabet the {\em Sturmian morphisms} on two letters. A morphism $\varphi$ is said to be {\em Sturmian} if $\varphi(\bs)$ is Sturmian for any Sturmian word $\bs$. The set of Sturmian morphisms over $\{a,b\}$ is closed under composition, and consequently it is a submonoid of the endomorphisms of $\{a,b\}^*$. Moreover, it is well known that the monoid of Sturmian morphisms is generated by the three morphisms: $(a\mapsto ab, b\mapsto a)$, $(a\mapsto ba. b \mapsto a)$, $(a\mapsto b, b \mapsto a)$ and that Sturmian morphisms are precisely the morphisms that map Sturmian words onto Sturmian words (see \cite{jBpS93acha, fMpS93morp}).

\subsection{Generators \& monoids}

By definition (see \cite{xDjJgP01epis, jJgP02epis}), the monoid of all {\em episturmian morphisms} $\cE$ is generated, under composition, by all the morphisms:  
\begin{itemize}
\item $\psi_a$:  $\psi_a(a) =a$, $\psi_a(x) =ax$ for any letter $x \neq a$;
\item $\bar \psi_a$:  $\bar \psi_a(a) =a$, $\bar \psi_a(x) =xa$ for any letter $x \neq a$;
\item $\theta_{ab}$: exchange of letters $a$ and $b$.
\end{itemize}
\begin{note} This system of generators is far from minimal, e.g., $\psi_a= \theta_{ab}\psi_b\theta_{ab}$, but gives simpler formulae.
\end{note} 

Moreover, the monoid  of so-called {\em epistandard morphisms} $\cS$ is generated by all the $\psi_a$ and the $\theta_{ab}$, and the monoid of {\em pure episturmian morphisms} $\cE_p$ (resp.~{\em pure epistandard morphisms} $\mathcal{S}_p$) is generated by the $\psi_a$ and $\bar{\psi}_a$ only (resp.~the $\psi_a$ only). The monoid $\cP$ of the {\em permutation morphisms} (i.e., the morphisms $\varphi$ such that $\varphi(\cA) = \cA$) is generated by all the $\theta_{ab}$. The importance of the monoid of pure episturmian morphisms will become clearer in the next section where we shall see that such morphisms are strongly linked to {\em spinned directive words} of episturmian words, which can be viewed as infinite sequences of rules for decomposing episturmian words by morphisms (see Theorems~\ref{T:epistandard-morph} and \ref{T:epi-morph}, to follow). In particular, any episturmian word is the image of another episturmian word by some pure episturmian morphism. 

The following diagram illustrates the inclusions between the monoids defined above. 

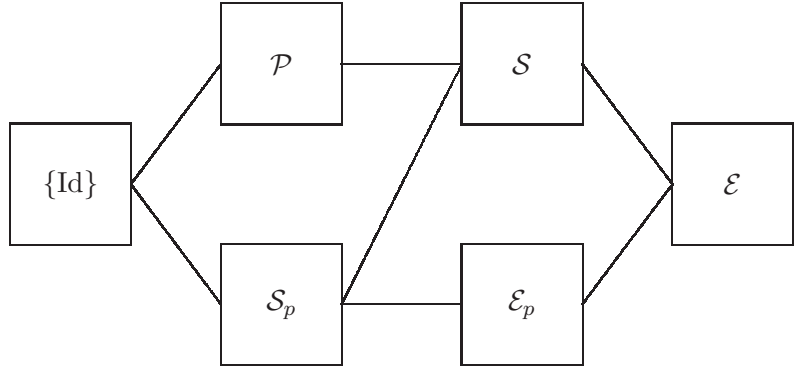
\unitlength = 0.8mm
\begin{center}
\begin{picture}(120,100)(15,0)

\put(10,50){\framebox(20,20){$\{\mathrm{Id}\}$}}
\put(45,30){\framebox(20,20){$\cS_p$}}
\put(45,70){\framebox(20,20){$\cP$}}
\put(85,30){\framebox(20,20){$\cE_p$}}
\put(85,70){\framebox(20,20){$\cS$}}
\put(120,50){\framebox(20,20){$\cE$}}
\qbezier(30,60)(30,60)(45,80)
\qbezier(65,80)(65,80)(85,80)
\qbezier(105,80)(105,80)(120,60)
\qbezier(30,60)(30,60)(45,40)
\qbezier(65,40)(65,40)(85,80)
\qbezier(65,40)(65,40)(85,40)
\qbezier(105,40)(105,40)(120,60)
\put(70,20){\makebox(0,0){Semidirect products: $ \cS= \cS_p \rtimes \cP, \quad \cE = \cE_p \rtimes \cP $}}
\end{picture}
\vspace{-1.2cm}
\end{center}

We note in particular that the monoid $\cE$ is a semidirect product of the submonoids of its pure morphisms and of its permutations. Consequently, any episturmian morphism $\varphi \in  \cE$ can be expressed in a unique way as $\varphi = \pi\mu = \mu'\pi$, where $\mu$, $\mu'$ are pure episturmian morphisms and $\pi$ is a permutation.

\begin{note} The episturmian morphisms are exactly the Sturmian morphisms when $|\cA| = 2$. 
\end{note}

Clearly, all episturmian morphisms on $\cA$ can be viewed as automorphisms of the free group generated by $\cA$ (e.g., see \cite{aG06onst, eG07repr, gR03conj, zWyZ99some}) and it follows that they are injective and that the monoids $\cE$ and $\cS$ are {\em left cancellative} (see \cite[Lemma~7.2]{gR03conj}) which means that for any episturmian morphisms $f, g, h$, if $fg = fh$ then $g = h$. Other fundamental properties of episturmian morphisms will be discussed in the next section and in Section~\ref{S:shifts&spins}. For an in-depth study of some further properties of these morphisms, the interested reader is referred to Richomme's paper \cite{gR03conj}, in which he considers invertibility, presentation, cancellativity, unitarity, characterization by conjugacy, and so on. Most of the results in \cite{gR03conj} naturally generalize those already known for Sturmian morphisms, but some new ones are also proved, such as a characterization of episturmian morphisms that preserve palindromes. In \cite{gR03lynd, gR04onmo}, Richomme also characterized the episturmian morphisms that preserve finite and infinite {\em Lyndon words} and those that preserve a lexicographic order on words.

\subsection{Relation with episturmian words} \label{SS:relation}

We now state two insightful characterizations of epistandard and episturmian words, which show that any episturmian word can be {\em infinitely decomposed} over the set of pure episturmian morphisms.


In the `standard' case:

\begin{theo} {\em \cite[Corollary 2.7]{jJgP02epis}} \label{T:epistandard-morph} An infinite word $\bs \in \cAw$ is \textbf{epistandard}  if and only if there exists an infinite word $\Delta = x_1x_2\cdots$ over $\cA$ and a sequence $(\bs ^{(i)})_{i \ge 0}$ of recurrent infinite words such that $\bs ^{(0)}=\bs$ and $\bs ^{(i-1)}=\psi_{x_i} (\bs^{(i)})$ for $i>0$.  
\qed
\end{theo}

In~\cite{jJgP02epis}, Justin and Pirillo showed that the infinite word $\Delta$ appearing in the above theorem is exactly the directive word of $\bs$ that arises from the equivalent definition of epistandard words given in Theorem~\ref{T:main}. In the binary case, the directive word $\Delta$ is related to the continued fraction expansion of the slope of the straight line represented by a standard word (see Chapter 2 in \cite{mL02alge}). 

\begin{exa} Recall the Tribonacci word $\br$, which has directive word $\Delta = (abc)^\omega$. We have 
$\br=\psi_a(\br^{(1)})$, where $\br^{(1)}$ is directed by $\TT(\Delta) = (bca)^\omega$. Notice that $\br^{(1)}=\pi(\br)$ with $\pi = (abc)$; a very particular case.
\end{exa}


More generally, the following result (Theorem~\ref{T:epi-morph}) extends the notion of a directive word to all episturmian words. Before stating the theorem, we need to introduce some more notation. First we define a new alphabet, $\bar{\cA} := \{\bar{x} ~|~ x \in \cA \}$. A letter $\bar x \in \bar\cA$ is considered to be $x$ with {\em spin} $1$,  whilst $x$ itself has {\em spin} $0$. The notion of a spin provides a convenient way to call upon the elementary pure episturmian morphisms $\psi_x$ and $\bar\psi_x$. Moreover, as well shall see in Section~\ref{S:shifts&spins}, it allows us to derive many properties of episturmian words from episturmian morphisms (as a consequence of the next theorem). This approach is used for instance in \cite{vBcHlZ06init, aG06acha, fLgR07quas,gR07aloc,gR07conj,rRlZ00agen} and of course in the papers of Justin {\em et al.}

A finite or infinite word over $\cA \cup \bar{\cA}$ is said to be a {\em spinned} word. Given a finite or infinite word $w = x_1x_2\cdots$ over $\cA$, we sometimes denote by $\breve w = \breve x_1 \breve x_2 \cdots$ any spinned word such that $\breve x_i = x_i$ if $x_i$ has spin $0$ and $\breve x_i = \bar x_i$ if $x_i$ has spin $1$. Such a word $\breve w$ is called a {\em spinned version} of~$w$. 

\begin{theo} {\em \cite[Theorem 3.10]{jJgP02epis}} \label{T:epi-morph}
An infinite word $\bt \in \cAw$ is \textbf{episturmian} if and only if there exists a spinned infinite word $\breve\Delta = \breve x_1 \breve x_2 \breve x_3 \cdots$ over $\cA \cup \bar\cA$ and an infinite sequence $(\bt^{(i)})_{i\geq0}$ of recurrent infinite words such that
\[
  \bt ^{(0)}=\bt \quad \mbox{and} \quad \bt ^{(i-1)} = \psi_{x_i} (\bt^{(i)}) \quad \mbox{or} \quad 
  \bt ^{(i-1)} =  \bar \psi_{x_i} (\bt^{(i)}) \quad \mbox{for all $i > 0$}, 
\]  
according to the spin $0$ or $1$ of $\breve x_i$, respectively.  
\end{theo}

For any epistandard word (resp.~episturmian word) $\bt$ and infinite word $\Delta$ (resp.~spinned infinite word $\breve\Delta$) satisfying the conditions of the Theorem~\ref{T:epistandard-morph} (resp.~Theorem~\ref{T:epi-morph}), we say that $\Delta$ (resp.~$\breve\Delta$) is a {\em directive word} (resp.~a {\em (spinned) directive word}) for $\bt$ or $\bt$ is {\em directed by} $\Delta$ (resp.~$\breve\Delta$). 

\begin{rem} \label{R:directive} It follows immediately from Theorem~\ref{T:epi-morph} that if $\bt$ is an episturmian word directed by a spinned infinite word $\breve\Delta$, then each $\bt^{(n)}$ (as defined in the theorem) is an episturmian word directed by $\TT^n(\breve\Delta) = \breve x_{n+1}\breve x_{n+2} \breve x_{n+3} \cdots$.
\end{rem}

The following important fact links Theorems~\ref{T:epistandard-morph} and \ref{T:epi-morph}. 

\begin{remq} \label{R:episturmian} \cite{jJgP02epis} If $\bt$ is an episturmian word directed by a spinned version $\breve \Delta$ of an infinite word $\Delta$ over $\cA$, then $\bt$ is factor-equivalent to the (unique) epistandard word $\bs$ directed by $\Delta$. 
\end{remq}

Moreover, with the same notation as in the above remark, the episturmian word $\bt$ is periodic if and only if the epistandard word $\bs$ is periodic, and this holds if and only if $|\Ult(\Delta)| = 1$  (see Remark~\ref{R:periodic} or Theorem~\ref{T:periodic} later). 

\begin{exa} Consider the episturmian word $\bm=baabacabab \cdots$ directed by $\breve \Delta =\bar a b \bar c (abc)^\omega$. Observe that $\bm$ is factor-equivalent to the Tribonacci word $\br$, and we have 
\[
\bm=\bar \psi_a(\bm^{(1)})= \bar \psi_a \psi_b (\bm^{(2)}) = \bar \psi_a \psi_b \bar \psi_c(\bm^{(3)}),
\]
where $\bm^{(3)}$ is directed by $\mathrm{T}^3 (\breve \Delta) = (abc)^\omega$, i.e., $\bm^{(3)} = \br$.
\end{exa}

\begin{exa} \label{ex:episkew} We now consider an example where the condition that the $\bt^{(i)}$ in Theorem~\ref{T:epi-morph} are recurrent is not satisfied. Let $\bt = d\br = dabacabaabacaba\cdots$ where $\br$ is the Tribonacci word and $d$ is a letter. Then $\bt = \bar\psi_a(\bt^{(1)})$, $\bt^{(1)} = \bar\psi_b(\bt^{(2)})$, $\bt^{(2)} = \bar\psi_c(\bt^{(3)})$, and so on; however, these $\bt^{(i)}$ are not recurrent (and $\bt$ is not episturmian). The infinite word $\bt = d\br$ is actually an example of an {\em episkew word}, i.e., a non-recurrent infinite word having episturmian factors. Such words are discussed in more detail in Section~\ref{SS:episkew}. 
\end{exa}

\begin{rem} Let us point out that the construction of epistandard words by palindromic closure (given in Theorem~\ref{T:main}) extends to all episturmian words: when $\breve x_n= \bar x_n$ write $x_n$ on the \emph{left} and use {\em palindromic left-closure}.  Here $\bm$ (from the above example) appears step by step on the right:
\begin{align*} 
 & \underline{a} \:\cdot \\
 & a  \cdot \underline{b}a \\
  aba\underline{c}&a \cdot ba \\
  abac&a  \cdot ba\underline{a}bacaba 
\end{align*}
\end{rem}

When an episturmian word is aperiodic, we have the following fundamental link between the words $(\bt^{(n)})_{n \geq 0}$ and the spinned infinite word $\breve\Delta$ occurring in Theorem~\ref{T:epi-morph}: if $a_n$ is the first letter of $\bt^{(n)}$, then $\mu_{\breve x_1 \cdots \breve x_n}(a_n)$ is a prefix of $\bt$ and the sequence $(\mu_{\breve x_1 \cdots \breve x_n}(a_n))_{n \geq 1}$ is not ultimately constant (since $\breve\Delta$ is not ultimately constant), then $\bt = \lim_{n \rightarrow \infty} \mu_{\breve x_1 \cdots \breve x_n}(a_n)$. This fact is a slight generalization of a result of Risley and Zamboni \cite[Prop.~III.7]{rRlZ00agen} on {\em S-adic representations} for standard Arnoux-Rauzy sequences.  See also the recent paper \cite{vBcHlZ06init} for S-adic representations of Sturmian words. Note that {\em $S$-adic dynamical systems} were introduced by Ferenczi \cite{sF99comp} as {\em minimal dynamical systems}  (e.g., see \cite{nP02subs})  generated by a finite number of substitutions. In the case of episturmian words, the notion itself is actually a reformulation of the well-known {\em Rauzy rules}, as studied in \cite{gR85mots}. In fact, it is well known that the {\em subshift} of an aperiodic episturmian word $\bt$ (i.e., the topological closure of the shift orbit of $\bt$)  is a {minimal dynamical system}, i.e., it  consists of all  the episturmian words with the same set of factors as $\bt$.

It is not hard to see that a morphism is episturmian (resp.~epistandard) if and only if it preserves the set of aperiodic episturmian (resp.~epistandard) words (see \cite{jJgP02epis}). Even more:

\begin{theo} {\em \cite[Theorem 3.13]{jJgP02epis}} \label{T:strict-morph}
A morphism $\varphi$ is episturmian (resp.~epistandard) if there exist strict episturmian (resp.~epistandard) words $\bm$, $\bt$ such that $\bm =\varphi(\bt)$. \qed 
\end{theo}

Purely morphic episturmian words (i.e., those generated by morphisms) are discussed further in Section \ref{S:shifts&spins}, where we consider the relationship between spins and the shifts that they induce. These ideas were used in \cite{jJgP04epis} to obtain a complete answer to the question: if an episturmian word is purely morphic, which shifts of it, if any, are also purely morphic? (See Theorem \ref{T:morphic}, to follow.) Such rigidity issues are discussed in more detail in Sections~\ref{SS:rigidity} and \ref{S:conclusion}. 

 In \cite{jJgP04epis}, Justin and Pirillo also made use of bi-infinite words, which often allow for more natural formulations. Indeed, the characterization (Theorem \ref{T:epi-morph}) of right-infinite episturmian words 
by a sequence $(\bt^{(i)})_{i\geq0}$ extends to bi-infinite episturmian words,  
with all the $\bt^{(i)}$ now bi-infinite episturmian words. That is, as for right-infinite episturmian words, we have bi-infinite words of the form $\bl^{(i)}.\br^{(i)}$ where $\bl^{(i)}$ is a left-infinite episturmian word and $\br^{(i)}$ is a right-infinite episturmian word. Moreover, if the bi-infinite episturmian word $\bb = \bl.\br$ is directed by $\breve\Delta$ with associated bi-infinite episturmian words $\bb^{(i)} = \bl^{(i)}.\br^{(i)}$, then $\br$ is directed by $\breve\Delta$ with associated right-infinite episturmian words $\br^{(i)}$.

\section{Spins, shifts, and directive words} \label{S:shifts&spins}

In this section, we discuss in more detail the notion of spins, the {\em shifts} they induce, and the concept of {\em block-equivalence} in connection with directive words. These notions allow us to study in particular when two different spinned infinite words direct the same episturmian word. Indeed, as we shall see in Section~\ref{SS:directive}, the correspondence between episturmian 
words and spinned directive words is not one-to-one.

\subsection{Notation for pure episturmian morphisms} \label{SS:pure}

For  $a \in \cA$, let $\mu_a = \psi_a$ and $\mu_{\bar a} = \bar\psi_a$. This operator $\mu$ can be naturally extended (as done in~\cite{jJgP02epis}) to a pure episturmian morphism: for any spinned finite word $\breve w = \breve x_1\cdots \breve x_n$ over $\cA \cup \bar\cA$, we define $\mu_{\breve w} := \mu_{\breve x_1}\cdots \mu_{\breve x_n}$ and set $\mu_\varepsilon$ equal to the identity morphism Id. 

Viewing $w = x_1x_2 \cdots x_n$ as a prefix of the directive word $\Delta=x_1x_2x_3\cdots \in \cAw$, it is clear from Theorem~\ref{T:epistandard-morph} that the words
\[
\mu_{x_1\cdots x_{n-1}}(x_{n}), \quad n \geq 1, 
\]
are prefixes of the epistandard word $\bs$ directed by $\Delta$. 

\begin{exa} We observe that any epistandard word $\bs \in \cAw$ has the form $\bs = \mu_w(\bs')$ for some uniquely determined finite word $w$ and strict epistandard word $\bs'$.  Indeed, if $\Delta = x_1x_2x_3\cdots \in \cAw$ is the directive word of $\bs$ and $m$ is the smallest positive integer such that $\Alph(x_{m+1}x_{m+2}\cdots) = \Alph(\Delta)$, then $x_1\cdots x_m$ is the shortest prefix of $\Delta$ that contains all the letters not appearing infinitely often in $\Delta$. Moreover, by Theorem~\ref{T:epistandard-morph}, $\bs = \mu_{x_1\cdots x_m}(\bs^{(m)})$ where $\bs^{(m)}$ is the epistandard word directed by $\TT^m(\Delta) = x_{m+1}x_{m+2}\cdots$. Since $\Ult(\TT^m(\Delta)) = \Alph(\TT^m(\Delta))$ by construction, the epistandard word $\bs^{(m)}$ is strict. For example, with $\Delta = c(ab)^\omega$, we have $\bs = \psi_c(\bs^{(1)})$ where $\bs^{(1)}$ is directed by $(ab)^\omega$, i.e., $\bs^{(1)}$ is the well-known \emph{Fibonacci word} over $\{a,b\}$. 
\end{exa}

For $n\geq 1$, let $u_{n+1} := Pal(x_1\cdots x_n)$ and set $u_1 = \empt$. Then by part $iv)$ of Theorem~\ref{T:main}, the epistandard word $\bs$ directed by $\Delta$ is given by $\bs = \lim_{n\rightarrow\infty} u_n$.  We have the following useful formula from \cite{jJgP02epis}:
\begin{equation} \label{eq:u_n**}
  u_{i+1} = \mu_{x_1\cdots x_{i-1}}(x_i) u_{i} \quad \mbox{for $i > 0$}.
\end{equation} 
For letters $(x_j)_{1 \leq j \leq i}$, formula~\eqref{eq:u_n**} inductively leads to:
\begin{equation}  \label{eq:u_n2}
u_{i+1} = \mu_{x_1 \cdots x_{i-1}}(x_i) \cdots \mu_{x_1}(x_2) x_1 = \prod_{1 \leq j \leq i}\mu_{x_1 \cdots x_{j-1}}(x_j).
\end{equation}
(Note that by convention, $x_1 \cdots x_0=\empt$ in the above product.) For example, with $\Delta = abcb\cdots$, we compute:
\[
  u_3 = Pal(abcb) = \mu_{abc}(b)\mu_{ab}(c)\mu_{a}(b)a = abacab\cdot abac \cdot ab\cdot a.
\]

\subsection{Shifts}

Now let $\breve w = \breve x_1 \breve x_2 \cdots \breve x_n$ be a spinned version of $w = x_1x_2\cdots x_n$ (viewed as a prefix of a spinned version $\breve \Delta$ of $\Delta$). Then, for any finite word $v$, we have
\begin{equation} \label{eq:f-shift}
  \mu_{\breve w}(v) = S_{\breve w}^{-1}\mu_w(v)S_{\breve w}  \quad \mbox{where $S_{\breve w} = \underset{\underset{\mid \breve x_i=\bar x_i}{i=n, \ldots, 1}}{\prod} \mu_{x_1 \cdots x_{i-1}}(x_i)$.}
\end{equation}

Observe that $S_{\breve w}$ is a prefix of $Pal(w)$; in particular $S_{\bar w} = Pal(w)$ by equation~\eqref{eq:u_n2}. Note also that $\mu_{\breve w}(v)= \TT^{|S_{\breve w}|}(\mu_w(v))$. The word $S_{\breve w}$ is called the {\em shifting factor} of $\mu_{\breve w}$ and its length $|S_{\breve w}|$ is called the {\em shift induced} by the prefix $\breve w$ of $\breve\Delta$ of length $n$ \cite{jJgP04epis}. 

\begin{exa}
If we take $\breve w=a \bar b c \bar a$, then 
\[
S_{\breve w}=\mu_{abc}(a) \mu_a(b)=abacaba\cdot ab. 
\]
Thus since $\mu_{abca}(ca) = abacabaab\cdot acabacaba$, we have 
\[
\mu_{a\bar b c \bar a}(ca) = \TT^9(\mu_{abca}(ca)) = acabacaba\cdot abacabaab.
\]   
\end{exa}

Likewise, for any infinite word $\by \in \cAw$, $\mu_{\breve w}(\by) = S_{\breve w}^{-1}\mu_w(\by)$. For example, if we take $\breve w = \bar a \bar b$, then $S_{\breve w} = Pal(ab) = aba$, and hence $\mu_{\bar a\bar b}(\by) = (aba)^{-1}\mu_{ab}(\by)$ for any infinite word $\by$.

\begin{note} The morphisms $\mu_w$ and $\mu_{\breve w}$ are {\em conjugate} morphisms \cite{gR03conj}.
\end{note}

\subsection{Block-equivalence \& directive words}  \label{SS:directive}

By Theorem~\ref{T:main} (and also Theorem~\ref{T:epistandard-morph}), any epistandard word $\bs \in \cAw$ has a unique directive word over $\cA$, but $\bs$ also has infinitely many other spinned directive words (see \cite{jJgP02epis,jJgP04epis, aGfLgR08dire}). For example, the Tribonacci word is  directed by $(abc)^\omega$ and also by $(abc)^n\bar{a}\bar{b}\bar{c}(a\bar{b}\bar{c})^\omega$ for each $n\geq 0$, as well as infinitely many other spinned words. The natural question: ``does any spinned word direct a unique episturmian word?'' was answered in \cite{jJgP02epis}.

\begin{prop}{\rm \cite{jJgP02epis}} 
\label{P:uniqueDirective}
\begin{enumerate}
\item  Any spinned infinite word $\breve \Delta$ having infinitely many letters with spin $0$ directs a unique episturmian word beginning with the left-most letter having spin $0$ in $\breve\Delta$. 
\item Any spinned infinite word $\breve \Delta$ with all spins ultimately $1$ directs exactly $|\Ult(\Delta)|$ episturmian words.
\item Let $\breve \Delta$ be a spinned infinite word having all its letters with spin $1$ and let $a \in \Ult(\Delta)$. Then $\breve \Delta$ directs exactly one episturmian word starting with $a$. \qed
\end{enumerate} 
\end{prop}

\begin{note} The above statement corrects a small error in Proposition~3.11 of~\cite{jJgP02epis} where item~3 was stated in the more general case when $\breve\Delta$ has all spins ultimately $1$. In this case, $\breve \Delta$ still directs exactly one episturmian word for each letter  $a$ in $\Ult(\Delta)$, but contrary to what is written in \cite{jJgP02epis}, nothing can be said about its first letter. 
\end{note}

{\em Block-equivalence} for spinned words was introduced in \cite{jJgP04epis} as a way of studying when $\breve \Delta$ and $ \hat \Delta$ (two spinned versions of a directive word $\Delta$) direct the same bi-infinite episturmian word. We do not recall the full details here, only a few notions relating to it.

\begin{notation} 
If $v \in \cA^+$, then $\bar v \in \bar\cA^+$ is $v$ with all spins $1$.
\end{notation}

A word of the form $xvx$, where $x \in \cA$ and $v \in (\cA\setminus\{x\})^*$, is called a ($x$-based) {\em block}. A ($x$-based) {\em block-transformation} is the replacement in a spinned word of an occurrence of $xv\bar x$ (where $xvx$ is a block) by $\bar x \bar v x$ or vice-versa. Two finite spinned words $w$, $w'$ are said to be {\em block-equivalent} if  we can pass from one to the other by a (possibly empty) chain of block-transformations, in which case we write $w \equiv w'$. For example, $\bar b\bar a b \bar c b \bar a \bar c$ and  $b a bc\bar b\bar a \bar c$ are block-equivalent because $\bar b\bar a b \bar c b \bar a \bar c \rightarrow ba\bar b \bar c b\bar a \bar c  \rightarrow b a bc\bar b\bar a \bar c$ and vice-versa. Note that if $w \equiv w'$ then $w$ and $w'$ are spinned versions of the same word over $\cA$. Block-equivalence extends to (right-)infinite words as follows.

Let $\Delta_1$, $\Delta_2$ be spinned versions of $\Delta$. We write $\Delta_1 \rightsquigarrow \Delta_2$ if there exist infinitely many prefixes ${f}_i$ of $\Delta_1$ and ${g}_i$ of $\Delta_2$  with the $g_i$ of strictly increasing lengths, and such that, for all $i$, $|g_i| \leq |f_i|$ and  ${f}_i \equiv {g}_i {c}_i$ for a suitable spinned word ${c}_i$. 
Infinite words $\Delta_1$ and $\Delta_2$
are said to be \textit{block-equivalent} (denoted by $\Delta_1 \equiv \Delta_2$) if $\Delta_1 \rightsquigarrow \Delta_2$ and $\Delta_2 \rightsquigarrow \Delta_1$.

\begin{rem} If $x$ is a letter and $v\in \cAstar $ is $x$-free, then $\bar x \bar v x$ and $x v \bar x$ are block-equivalent and they induce the same shift, i.e., $\mu _{\bar x \bar v x} = \mu _{ x v \bar x}$ \cite[Theorem 2.2]{jJgP04epis}. Thus the monoid of pure episturmian morphisms, $\cE_p$, is isomorphic to the quotient of $(\cA \cup \bar\cA)^*$ by the block-equivalence generated by
\[ 
\{\bar x \bar v x \equiv xv \bar x  \mid x \in \cA, v \ \mathrm{is} \  \mbox{$x$-free} \}.
\] 
Note that this has some relation to the study of conjugacy and episturmian morphisms carried out by Richomme \cite{gR03conj}. 
\end{rem}

From what we have already learned about bi-infinite episturmian words (in Sections~\ref{SS:bi-infinite} and \ref{SS:relation}), it is clear that Justin and Pirillo's results about spinned infinite words directing the same bi-infinite episturmian word are still valid for words directing the same (right-infinite) episturmian word. Roughly speaking, two spinned infinite words direct the same episturmian word if and only if they are block-equivalent. For instance, we have the following results for {\em wavy} spinned versions of $\Delta \in \cAw$. A spinned version $\breve \Delta$ of $\Delta$ is said to be {\em wavy} if $\breve \Delta$ contains infinitely many letters of spin~$0$ and infinitely many letters of spin~$1$. 

\begin{theo} {\em \cite[Theorem~3.4]{jJgP04epis}}
Suppose $\breve\Delta$ and $\hat\Delta$ are wavy versions of $\Delta \in \cAw$ with $|\Ult(\Delta)| > 1$. Then $\breve\Delta$ and $\hat\Delta$ direct the same episturmian word if and only if $\breve\Delta \equiv \hat\Delta$. \qed
\end{theo}

For example, $ba(\bar b c \bar a)^\omega$ and $\bar b \bar a b (c\bar  a \bar b)^\omega$ direct the same episturmian word, namely $\mu_{ba\bar b c}(c\br)$ ($= \mu_{\bar b \bar a b c}(c\br)$) where $\br$ is the Tribonacci word.

\begin{theo} {\em \cite[Prop.~3.6]{jJgP04epis}} \label{P:wavy-2} Let $\breve\Delta$, $\hat\Delta$ be two spinned versions of $\Delta \in \cAw$ with $|\Ult(\Delta)| > 1$, $\breve\Delta$ wavy, and $\hat\Delta$ having all spins ultimately $0$ or $1$. If $\breve\Delta$ and $\hat\Delta$ direct the same episturmian word, then $\breve\Delta \rightsquigarrow \hat\Delta$. \qed
\end{theo}

Similar results also hold when all spins are ultimately $0$ or $1$ and in the periodic case. See Propositions 3.7 and 3.10 in \cite{jJgP04epis}. 

\begin{rem} \label{R:ostrowski}
In \cite{jJgP04epis}, the study of block-equivalence for finite spinned words led to numeration systems that resemble the {\em Ostrowski systems} \cite{vB01auto} associated with Sturmian words. A matrix formula for computing the number of representations of an integer in such a system was also given in \cite[Section~2]{jJgP04epis}.
\end{rem}

 More recently, Glen, Lev\'e, and Richomme~\cite{aGfLgR08dire} established the following complete characterization of pairs of spinned infinite words directing the same unique episturmian word. Not only does the following theorem provide the relative forms of two spinned infinite words directing the same episturmian word, but it also fully solves the periodic case, which was only partially solved in~\cite{jJgP04epis}.

 \begin{theo}\label{T:directSame} {\em \cite{aGfLgR08dire}}
Given two spinned infinite words $\Delta_1$ and $\Delta_2$, the following assertions are equivalent.
\begin{description}
\item{i)} $\Delta_1$ and $\Delta_2$ direct the same right-infinite episturmian word.
\item{ii)} One of the following cases holds for some $i, j$ such that $\{i, j\} = \{1, 2\}$:
\begin{enumerate}
\item \label{Ti:1} $\Delta_i = \prod_{n \geq 1} v_n$, $\Delta_j = \prod_{n \geq 1} z_n$ where $(v_n)_{n \geq 1}, (z_n)_{n \geq 1}$ are spinned words such that $\mu_{v_n} = \mu_{z_n}$ for all $n \geq 1$;

\item \label{Ti:2} $\Delta_i = {w} x \prod_{n \geq 1} v_n {\breve x}_n$, $\Delta_j = {w'} {\bar x} \prod_{n \geq 1} {\bar v}_n {\hat x}_n$ where ${w}$, ${w'}$ are spinned words such that $\mu_{w} = \mu_{w'}$, $x$ is a letter, $(v_n)_{n \geq 1}$ is a sequence of non-empty $x$-free words, and $({\breve x}_n)_{n \geq 1}$, $({\hat x}_n)_{n \geq 1}$ are sequences of non-empty spinned words over $\{x, \bar x\}$ such that, for all $n \geq 1$, $|{\breve x}_n| = |{\hat x}_n|$ and  $|{\breve x}_n|_x = |{\hat x}_n|_x$;

\item \label{Ti:4} $\Delta_1 = w \bx$ and $\Delta_2 = w'\by$ where $w$, $w'$ are spinned words, $x$ and $y$ are letters,  and  $\bx \in \{x, \bar x\}^\omega$, $\by \in \{y, \bar y\}^\omega$ are spinned infinite words such that $\mu_{w}(x) = \mu_{w'}(y)$. 
\end{enumerate} \qed
\end{description} 
 \end{theo}
 
 In items \ref{Ti:1} and \ref{Ti:2} of Theorem~\ref{T:directSame}, the two considered directive words are spinned versions of the same infinite word. This does not hold in item \ref{Ti:4}, which concerns only periodic episturmian words. In particular, we observe the following:

\begin{rem} 
\label{R:aperiodic-spinned-versions}  If an {\em aperiodic} episturmian word is directed by two different spinned infinite words $\Delta_1$ and $\Delta_2$, then $\Delta_1$ and $\Delta_2$ are spinned versions of the same word $\Delta$. 
\end{rem}

As an example of item \ref{Ti:4}, one can consider the periodic episturmian word $(bcba)^\omega$ which is directed by both $bca^\omega$ and $b \bar ac^\omega$ (since $\mu_{bc}(a) = \mu_{b\bar a}(c)$). Note also that $(bcba)^\omega$ is epistandard and has the same set of factors as the epistandard word $(babc)^\omega$ directed by $bac^\omega$. Actually, in view of Remark~\ref{R:episturmian}, we observe the following:

\begin{rem} \label{R:subshift}  The subshift of any aperiodic episturmian word contains a unique (aperiodic) epistandard word, whereas the subshift of a periodic episturmian word contains exactly two (periodic) epistandard words, except if this word is $a^\omega$ with $a$ a letter. 
\end{rem}

We also observe that $x$ and $y$ can be equal in item~\ref{Ti:4} of Theorem~\ref{T:directSame}; for example $(ab)^\omega$ is directed by $a\bar b b^\omega$ and by $ab^\omega$.

\begin{exa} \cite{aGfLgR08dire}
For $a, b, c$ three different letters in $\cA$, the spinned infinite words $\Delta_1 = a(bc{\bar a})^\omega$ and $\Delta_2 = {\bar a}({\bar b}{\bar c}{\bar a})^\omega$ direct the same episturmian word that starts  with the letter $a$. Indeed, these two directive words fulfill item~2 of Theorem~\ref{T:directSame} with $w = w' = \varepsilon$, $x = a$, and for all $n$, $v_n = bc$ and $\breve{x}_n = \hat{x}_n = \bar{a}$. Moreover the fact that $\Delta_1$ starts with the letter $a$ shows that the word it directs starts with $a$. 
Similarly $\Delta_1' = {\bar a}b(ca{\bar b})^\omega$ and $\Delta_2' = {\bar a}{\bar b}({\bar c}{\bar a}{\bar b})^\omega$ direct the same episturmian word starting with the letter $b$.
Since $\Delta_2 = \Delta_2'$,   
this shows that the relation ``direct the same episturmian word'' over spinned infinite words is not an equivalence relation. 
\end{exa}

Items \ref{Ti:2} and \ref{Ti:4} of Theorem~\ref{T:directSame} show that any episturmian word is directed by a spinned infinite word having infinitely many letters of spin~0, but also by a spinned word having both infinitely many letters of spin~0 and infinitely many letters of spin~1 (i.e., a wavy word). To emphasize the importance of these facts, let us recall from Proposition~\ref{P:uniqueDirective} that if $\breve\Delta$ is a spinned infinite word over $\cA \cup \bar\cA$ with infinitely many letters of spin~0, then there exists a unique episturmian word $\bt$ directed by $\breve\Delta$. Unicity comes from the fact that the first letter of $\bt$ is fixed by the first letter of spin~0 in $\breve\Delta$. We also note that if an episturmian word $\bt$ has two directive words satisfying items \ref{Ti:2} or \ref{Ti:4}, then $\bt$ has infinitely many directive words (this was shown in \cite{aGfLgR08dire}).

When studying repetitions in Sturmian words, Berth\'e, Holton, and Zamboni~\cite{vBcHlZ06init} proved that any Sturmian word has a unique directive word over $\{a,b,\bar a, \bar b\}$ containing infinitely many letters of spin~0, but no factor of the form ${\bar a}{\bar b}^na$ or ${\bar b}{\bar a}^nb$ with $n$ an integer. Lev\'e and Richomme~\cite{fLgR07quasB} recently generalized this result  to episturmian words by introducing a way to `normalize' the directive word(s) of an episturmian word so that any episturmian word can be defined uniquely by its so-called {\em normalized directive word}, defined by some factor avoidance, as follows. This idea has since proved useful in the study of {\em quasiperiodic} episturmian words (see Section~\ref{S:conclusion}); in particular, it provides an effective way to decide whether or not a given episturmian word is quasiperiodic.

\begin{theo}\label{T:normalisation} {\em \cite{aGfLgR08dire, fLgR07quasB}}
Any episturmian word $\bt \in \cAw$ has a spinned directive word $\breve\Delta$ containing infinitely many letters of spin $0$, but no factor in ${\bigcup}_{a\in\cA} \bar a \bar\cA^*a$. Such a directive word is unique if $\bt$ is aperiodic, in which case $\breve\Delta$ is called the {\em normalized directive word} for $\bt$. \qed
\end{theo}
\begin{note} Unicity does not necessarily hold for periodic episturmian words. For example, the periodic  episturmian word $(ab)^\omega = \psi_a(b^\omega) = \bar\psi_b(a^\omega)$ is directed by $a{b}^\omega$ and by $\bar{b}a^\omega$ (since $\psi_a(b) = ab = \bar\psi_b(a)$). 
\end{note}

The following result tells us precisely which episturmian words have a unique directive word.

\begin{theo} \label{T:uniqueDirective} {\em \cite{aGfLgR08dire}}
An episturmian word $\bt \in \cAw$ has a unique directive word if and only if the (normalized) directive word of $\bt$ contains 1) infinitely many letters of spin $0$, 2) infinitely many letters of spin $1$, 3) no factor in ${\bigcup}_{a\in\cA} \bar a \bar\cA^*a$, and 4) no factor in ${\bigcup}_{a\in\cA} a \cA^*\bar a$. Such an episturmian word is necessarily aperiodic. \qed
\end{theo}

For instance, a particular family of episturmian words having unique directive words consists of those directed by {\em regular wavy words} \cite{aG07orde,aGfLgR07quas}, i.e., spinned infinite words having both infinitely many letters of spin $0$ and infinitely many letters of spin $1$ such that each letter occurs with the same spin everywhere in the directive word.  More formally, a spinned version $\breve w$ of a finite or infinite word $w$ is said to be {\em regular} if, for each letter $x \in \Alph(w)$, all occurrences of $\breve x$ in $\breve w$ have the same spin $(0$ or $1)$. For example, the regular wavy word $(a\bar b \bar c)^\omega$  is the unique directive word for the episturmian word $a\br =  a{a}{b}a{c}aba{a}bacab\cdots~$ where $\br$ is the Tribonacci word.

In the Sturmian case, we have:

\begin{prop} \label{P:uniqueDirective-Sturmian} {\em \cite{aGfLgR08dire}}
Any Sturmian word has either a unique spinned directive word or infinitely many spinned directive words. Moreover, a Sturmian word has a unique directive word if and only if its (normalized) directive word is regular wavy. \qed
\end{prop}

As pointed out in \cite{aGfLgR08dire}, Proposition~\ref{P:uniqueDirective-Sturmian} shows a great difference between Sturmian words and episturmian words constructed over alphabets with at least three letters. Indeed, when considering words over a ternary alphabet, one can find episturmian words having exactly $m$ directive words  for any $m \geq 1$. For instance, the episturmian word $\bt$ directed by $\breve\Delta = a(b\bar a)^{m-1}b\bar c(ab\bar c)^\omega$ has exactly $m$ directive words, namely $(\bar a\bar b)^ia(b\bar a)^{j}b\bar c(ab\bar c)^\omega$ with $i+j = m-1$. Notice that the suffix $b\bar c(ab\bar c)^\omega$ of $\breve\Delta$ is regular wavy, and the other $m-1$ spinned versions of $\Delta$ that also direct $\bt$ arise from the $m-1$ words that are block-equivalent to the prefix $a(b\bar a)^{m-1}$.

\subsection{Periodic and purely morphic episturmian words} \label{SS:rigidity}

We are now ready to describe periodic and purely morphic episturmian words.

Recall from Remark \ref{R:periodic} that the periodic episturmian words correspond to $|\Ult(\Delta)| = 1$. The following theorem gives the form of such words in terms of pure episturmian morphisms.

\begin{theo} {\em \cite{jJgP02epis}}  \label{T:periodic}  An episturmian word is periodic if and only if it is $(\mu_{\breve w}(x))^\omega$ for some spinned finite word $\breve w$ and letter $x$. \qed
\end{theo}

For example, $(\mu_{a\bar{b}}(c))^\omega  = (acab)^\omega$ is the periodic episturmian word directed by $a\bar bc^\omega$ (in fact, it is epistandard as it is also directed by $acb^\omega$).


The next theorem characterizes purely morphic episturmian words with respect to their directive words.

\begin{theo} {\em \cite[Theorem~3.14]{jJgP02epis}} \label{T:purely-morphic}
An aperiodic episturmian word is purely morphic (i.e., generated by a morphism) if and only if it is directed by a periodic spinned infinite word $\breve \Delta = (\breve f)^\omega$ for some spinned word $\breve f$. Moreover it can be  generated by~$\mu_{\breve f}$. \qed
\end{theo}
We observe from Theorem~\ref{T:purely-morphic} that any purely morphic episturmian word is strict (i.e., an Arnoux-Rauzy sequence) as $\Ult(\Delta) = \Alph(f) = \Alph(\Delta)$. The proof of this theorem makes use of Proposition~\ref{P:uniqueDirective} and Theorem~\ref{T:strict-morph}.

\begin{exa} The Tribonacci word is generated by $\mu_{abc}$. Notice that $\mu_{abc} = \sigma^3$ where $\sigma$ is the {\em Tribonacci morphism} defined by $\sigma: (a,b,c)\mapsto (ab,ac,a)$.
\end{exa}

\begin{rem} Purely morphic standard Sturmian words were previously characterized independently in the following papers: \cite{jBpS93acha, dCwMaPpS93subs, tKaV96subs}. Yasutomi~\cite{sY99onst}  has since established a characterization of all purely morphic Sturmian words with respect to their  slopes and intercepts (when viewed as cutting sequences). An alternative geometric proof of Yasutomi's result was recently given by Berth\'e {\em et al.} in \cite{vBhEsIhR07onsu}.
\end{rem}

Using the notion of block-equivalence, Justin and Pirillo \cite{jJgP04epis} explicitly determined which shifts, if any, of a purely morphic episturmian word are also purely morphic.

\begin{theo} {\em \cite{jJgP04epis}} \label{T:morphic}
If an episturmian word $\bt$ is purely morphic, then its shift $\TT^i(\bt)$ is also purely morphic if and only if $i$ belongs to some particular interval. \qed
\end{theo}
See Section 4 of \cite{jJgP04epis} for specific (and very technical) details. 

\begin{exa} For the Tribonacci word $\br$, only itself and $\TT^{-1}(\br)$ are purely morphic. Note that $\TT^{-1}(\br)$ corresponds to three episturmian words: $a\br$, $b\br$, $c\br$, directed by $(a \bar b \bar c)^\omega$, $(\bar a b \bar c)^\omega$, $(\bar a \bar b c)^\omega$, respectively.
\end{exa}

\begin{remq}
Theorem~\ref{T:morphic} corrects an error in \cite[Section~5.1]{jJgP02epis} where it was mistakenly said  that if an episturmian word is purely morphic then any shift of it is also purely morphic. Indeed, this is false even in the Sturmian case as Fagnot \cite{iF06alit}  has shown that if $\bs$ is a purely morphic standard Sturmian word on $\{a,b\}$, then $a\bs$, $b\bs$, $ab\bs$, $ba\bs$ (which are purely morphic \cite{jBpS94arem}) are the only purely morphic Sturmian words related to $\bs$ by a shift.
\end{remq}

\section{Arnoux-Rauzy sequences} \label{S:A-R}

We now briefly turn our attention to Arnoux-Rauzy sequences  since their combinatorial properties are also considered in the sections that follow.

{\em Arnoux-Rauzy sequences}  
are uniformly recurrent infinite words over a finite alphabet $\cA$ with factor complexity $(|\cA|-1)n + 1$ for each $n \in \NN$, and exactly one right and one left special factor of each length. They were introduced by Arnoux and Rauzy \cite{gR82suit, pAgR91repr}, who studied them using {\em Rauzy graphs}, 
with particular emphasis on the case $|\cA| = 3$. (Note that the foregoing definition is equivalent to the one given in the introduction.)

As mentioned previously (in Section~\ref{SS:strict}), Arnoux-Rauzy sequences are exactly the strict episturmian words; in particular, any episturmian word has the form $\varphi(\bt)$ with $\varphi$ an episturmian morphism and $\bt$ an Arnoux-Rauzy sequence. In this sense, episturmian words are only a slight generalization of Arnoux-Rauzy sequences. For example, the family of episturmian words on three letters $\{a,b,c\}$ consists of the Arnoux-Rauzy sequences over $\{a,b,c\}$, the Sturmian words over $\{a,b\}$, $\{b,c\}$, $\{a,c\}$ and their images under episturmian morphisms on $\{a,b,c\}$, and periodic infinite words of the form $\varphi(x)^\omega$ where $\varphi$ is an episturmian morphism on $\{a,b,c\}$ and $x \in \{a,b,c\}$.

Arnoux-Rauzy sequences have deep properties studied in the framework of dynamical systems, with connections to geometrical realizations such as {\em Rauzy fractals} \cite{pAsI01piso} and interval exchanges. When $|\cA| = 3$, the condition on the special factors distinguishes Arnoux-Rauzy sequences from other infinite words of complexity $2n+1$, such as those obtained by coding trajectories of $3$-interval exchange transformations (e.g., see \cite{sFcHlZ01stru}). In \cite{pAgR91repr}, it was shown how Arnoux-Rauzy sequences of complexity $2n+1$ (i.e., the $3$-strict episturmian words) can be geometrically realized by an exchange of six intervals on the unit circle, which generalizes the representation of Sturmian sequences by rotations.

An alternative way of introducing and studying Arnoux-Rauzy sequences is in the context of {\em $S$-adic dynamical systems}, as done in \cite{rRlZ00agen} for instance (see our remarks following Theorem~\ref{T:epi-morph} in Section~\ref{SS:relation}). In \cite{dDlZ03comb}, Damanik and Zamboni give a kind of survey on this approach by considering {\em Arnoux-Rauzy subshifts} and answering various combinatorial questions concerning linear recurrence, maximal powers of factors, and the number of palindromes of a given length. They also present some applications of their results to the spectral theory of discrete one-dimensional Schr$\ddot{\textrm{o}}$dinger operators with potentials given by Arnoux-Rauzy sequences.

Arnoux-Rauzy sequences also have interesting arithmetical properties. For instance, if one considers the {\em frequencies} of letters (as discussed later in Section~\ref{SS:frequencies}), they are well-defined, and renormalization by an episturmian morphism leads to a generalization of the continued fraction algorithm that associates to each $k$-letter Arnoux-Rauzy sequence an infinite array of $k \times k$ rational numbers.  In the special case $k = 2$, these fractions are consecutive Farey numbers arising from the continued fraction expansion of the frequencies of the two letters. More generally, given an Arnoux-Rauzy sequence on $k$-letters, its directive word is determined by the `multi-dimensional' continued fraction expansion of the frequencies of the first $k-1$ letters. Unfortunately, this generalized algorithm (except for the case $k=2$ when it is exactly the usual continued fraction algorithm) is only defined on a set of measure zero in $\RR^{k-1}$.  
This reduces its interest and explains why it has not been appropriately studied since its inception (see Sections~\ref{SS:palindromic-closure} and \ref{SS:frequencies} for further details). Nonetheless, a nice arithmetical characterization of $3$-letter Arnoux-Rauzy sequences can be given, as follows. We say that a triple $(a,b,c)$ does not satisfy the triangular   inequality if one of the coordinates is larger than the sum of the other two (e.g., $a > b+c$). ÊIn that case, we can renormalize in a unique way to obtain the triple $(a-b-c, b, c)$ satisfying the triangular inequality. The set of allowable frequencies for $3$-letter Arnouxy-Rauzy sequences is exactly the set of triples $(a,b,c)$ that can be infinitely renormalized, each time to a triple that does {\it not} satisfy the triangular inequality (see \cite{pAgR91repr}). The resulting picture exhibits a kind of {\em Sierpinski carpet}. 

For further details on Arnoux-Rauzy sequences, we refer the reader to the interesting survey \cite{vBsFlZ05inte} in which Berth\'e, Ferenczi, and Zamboni discuss connections between Arnoux-Rauzy sequences and rotations of the $2$-torus; coding of two-dimensional actions and two-dimensional Sturmian words; and interval exchanges and sequences of low complexity. See also  \cite{nCpHaM01prop}, Section~12.2.3 in \cite{nP02subs}, and J.~Berstel's nice survey paper \cite{jB02rece} in which he compares some combinatorial properties of Arnoux-Rauzy sequences (as well as  episturmian words) to those of Sturmian words. 


\subsection{Finite Arnoux-Rauzy words}

A finite word $w$ is said to be \emph{finite episturmian} if $w$ is a factor of some infinite episturmian  word. When considering factors of (infinite) episturmian words, it suffices to consider only the strict standard ones (i.e., the standard Arnoux-Rauzy sequences). Indeed, for any prefix $u$ of an epistandard word, there exists a strict epistandard word also having $u$ as a prefix.  
In particular, the words $\mu_w (x)$, with $w \in \cA^*$ and $x \in \cA$, are the {\it standard} ones ({\it cf.} standard words, e.g., \cite[Chapter 2]{mL02alge}). They can be obtained by the {\it Rauzy rules} \cite{gR85mots} (see also \cite[Theorem 8]{xDjJgP01epis}), and this has a strong connection with the set of periods of the palindromes $u_{n+1} = Pal(x_1\cdots x_n)$ (given in Theorem~\ref{T:main})  and the Euclidean  algorithm. This relation was studied by Castelli, Mignosi, and Restivo \cite{mCfMaR99fine}, who extended the well-known {\em Fine and Wilf Theorem} \cite{mL83comb} 
to words having three periods. Justin~\cite{jJ00onap} generalized this result even further to words having an arbitrary number of periods, which led  to a characterization of finite episturmian words.

Finite episturmian words are exactly the {\em finite Arnoux-Rauzy words}. Such words were enumerated by Mignosi and Zamboni \cite{fMlZ02onth}, who described a multi-dimensional generalization of the {\em Euler phi-function} that counts the number of finite Arnoux-Rauzy words of each length. Finite episturmian words have also been characterized with respect to lexicographic orderings in \cite{aGjJgP06char} (see Theorem \ref{T:aGjJgP06char} later).

\section{Some properties of factors} \label{S:factors}

\subsection{Factor complexity}

As mentioned previously, any $k$-strict episturmian word has complexity $(k-1)n+1$ for all $n \in \NN$. More generally:

\begin{theo} {\em \cite[Theorem~7]{xDjJgP01epis}}
Suppose $\bt$ is an episturmian word directed by $\breve\Delta$ with $|\Ult(\Delta)| > 1$. Then, for $n$ large enough, $\bt$ has complexity $(k - 1)n + q$ for some $q \in \NN^+$, where $k = |\Ult(\Delta)|$. \qed
\end{theo}
This theorem can be easily deduced from the fact that for sufficiently large $n$, any left special factor of  $\bt$ of length at least $n$ has exactly  $k = |\Ult(\Delta)|$ different left extensions in $\bt$ (by Theorem~6 in \cite{xDjJgP01epis}).

\subsection{Palindromic factors}

The {\em palindromic complexity} of episturmian words was established in \cite{jJgP02epis} by carrying out a similar study to the one for Sturmian words in \cite{xDgP99pali}. 

\begin{theo} {\em \cite[Theorem~4.4]{jJgP02epis}} If $\bt$ is an $\cA$-strict episturmian word, then there exists exactly 
 \begin{itemize}
 \item  one palindrome of length $n$ for all even $n$,  
 \item  one palindrome of length $n$ and centre $x$ for all odd $n$ and $x \in \cA$. \qed
\end{itemize}
\end{theo}
As shown in \cite{xDgP99pali}, the above property is characteristic in the Sturmian case, but not when $\cA$ contains more than two letters because it also holds for {\em billiard words}, which are not episturmian (see Borel and Reutenauer \cite{jBcR05pali}).

\begin{theo} {\em \cite[Section~4.2]{jJgP02epis}}
If $\bt$ is episturmian, then there exist $|\Ult(\Delta)|+1$ bi-infinite episturmian words of the form $\tilde\bm.\bm$ and $\tilde\bm x \bm$ with $x \in \Ult(\Delta)$ giving the palindromic factors of $\bt$. The spinned versions of $\Delta$ directing these bi-infinite episturmian words can be easily constructed via a simple algorithm.  \qed
\end{theo}
For more precise technical details, see Section 4.2 in \cite{jJgP02epis}.

\begin{exa}
For the Tribonacci word, $\tilde\br.\br$ is directed by $(abc \bar a bca \bar b cab \bar c)^\omega$.
\end{exa}

\subsubsection{Iterated palindromic closure} \label{SS:palindromic-closure}

In \cite{rRlZ00agen}, Risley and Zamboni gave an alternative construction of the sequence $(u_n)_{n\geq 1}$ of palindromic prefixes of an epistandard word (where $u_1 = \empt$ and $u_{i+1} = Pal(x_1\cdots x_{i})$ for all $i \geq 1$), using a `hat operation' as opposed to palindromic closure. The so-called  hat operation is defined as follows. We construct a new alphabet $\cA' := \cA \cup \widehat{\cA}$ where $\widehat{\cA} = \{\widehat x ~|~ x \in \cA\}$ and denote by $\phi$ the morphism $\phi: \cA' \to \cA$ defined by $\phi(x) = \phi(\widehat x) = x$ for all letters $x \in \cA$. The morphism $\phi$ extends to a morphism (also denoted by $\phi$) from words over $\cA'$ to words over $\cA$. Now, from a given directive word $\Delta = x_1x_2x_3\cdots \in \cAw$, we construct a sequence of words $(p_i)_{i\geq 1}$ as follows. We begin with $p_1 = \empt$ and $p_2 = \widehat x_1$. Then, for $n \geq 2$, $p_{n+1}$ is obtained from $p_{n}$ according to the rule: if $\widehat x_n$ does not occur in $p_{n}$, then $p_{n+1} = p_n \widehat x_n \phi(p_n)$; otherwise $p_{n+1} = p_n\widehat x_n\phi(s_n)$, where $s_n$ is the longest palindromic suffix of $p_n$ containing no occurrence of $\widehat x_n$. 

\begin{exa} Let $\Delta = (abc)^\omega$. Then using the hat operation, we obtain:
\begin{eqnarray*}
p_1 &=& \empt \\
p_2 &=& \hat a \\
p_3 &=& \hat a \hat b a \\
p_4 &=& \hat a \hat b a \hat c aba \\
p_5 &=& \hat a \hat b a \hat c aba \hat a bacaba \\
p_6 &=& \hat a \hat b a \hat c aba \hat a bacaba \hat bacabaabacaba \\
~ &\vdots& ~
\end{eqnarray*}
Now removing all hats (by applying $\phi$), we see that the $p_i$'s are precisely the palindromic prefixes of the Tribonacci word:  $abacabaabacababacabaabacaba\cdots$.
\end{exa}

As demonstrated by the above example, the hat operation is clearly the same as iterated palindromic closure; in fact, the relationship between these two constructions is evident by formula~\eqref{eq:u_n**}, which we now rewrite as:
\[
  Pal(x_1\cdots x_n) = \mu_{x_1\cdots x_{n-1}}(x_n)Pal(x_1\cdots x_{n-1}) \quad \mbox{for $n > 0$}. 
\]
The above formula is actually a special case of formula (3) from \cite{jJ05epis}, which also happens to be formula (3) in \cite{jJgP02epis}, namely:
\begin{equation} \label{eq:formula(3)}
Pal(vw) = \mu_v(Pal (w))Pal(v) \quad \mbox{for any words $w$, $v$}.
\end{equation}
This formula is commonly referred to as {\em Justin's Formula}, from which we deduce the following two  special cases:
\begin{equation} \label{eq:formula(3)*}
Pal(xw) = \psi_x(Pal(w))x \quad \mbox{and} \quad Pal(wx)=\mu_w(x)Pal(w) \quad \mbox{for any word $v$ and letter $x$. }
\end{equation}
The first formula given in \eqref{eq:formula(3)*} tells us that $Pal(xw)$ is obtained from $Pal(w)$ simply by inserting the letter $x$ before each letter different from $x$ and then appending $x$ to the resulting word. For example, $Pal(bc) = bcb$ and $Pal(abc) = abacaba$. The second formula given in \eqref{eq:formula(3)*} provides another way to compute the palindromic right-closure of $wx$ by placing the finite epistandard word $\mu_w(x)$ in front of $Pal(w)$. For example, to compute $Pal(abcb)$ we need only compute the words $\mu_{abc}(b) = abacab$ and $Pal(abc) = abacaba$, and then we have:
\[
  Pal(abcb) = \mu_{abc}(a)Pal(abc) = abacab\cdot abacaba.
\]

In \cite{jJ05epis}, Justin established some relations between the words $Pal(w)$, $\mu_w$, $Pal(\rev w)$, and $\mu_{\rev w}$ where $w$ is any finite word. Moreover, he showed that his results can be explained by the similarity of the {\em incidence matrices} of $\mu_w$ and $\mu_{\rev w}$. One curious result is that $|Pal(w)|=|Pal(\tilde w)|$. For example, with $w = abac$, $Pal(w) = abaabacabaaba$ and $Pal(\rev w) = cacbcacacacbcac$, both of length $15$.

Applying his results to a 2-letter alphabet, Justin \cite{jJ05epis} gave a new proof of a Galois theorem on continued fractions, by considering the epistandard words that are fixed points of $\mu_{w}$ and $\mu_{\rev w}$ for any finite word $w$.  From this point of view, Justin's result highlights the relevance of the previously  mentioned `multi-dimensional' continued fraction algorithm, proposed by Zamboni \cite{lZ98uneg, nWlZ01freq} (see also \cite[Section~12.2]{nP02subs}). However, there still remains much work to be done in this direction, especially concerning the {\em generalized intercept} (coherent with the Sturmian case)  introduced in \cite[Section~5.4]{jJgP02epis} and the generalized Ostrowski numeration systems \cite{vB01auto, jJgP04epis} (recall Remark \ref{R:ostrowski}).  \smallskip

\begin{note} The aforementioned Galois theorem was used in the theory of Sturmian words to characterize so-called {\em Sturm numbers} (see \cite[Theorem~2.3.26]{mL02alge}).
\end{note}

\subsubsection{Palindromic richness} \label{SS:rich}

In \cite{xDjJgP01epis}, Droubay, Justin, and Pirillo observed that any finite word $w$ contains at most $|w| + 1$ distinct palindromes (including the empty word). Even further, they proved that a word $w$ contains exactly $|w| + 1$ distinct palindromes if and only if the longest palindromic suffix of any prefix $p$ of $w$ occurs exactly once in $p$ (i.e., every prefix of $w$ has {\em Property $Ju$} \cite{xDjJgP01epis}). Such words  are `rich' in palindromes in the sense that they contain the maximum number of different palindromic factors. Accordingly, we say that a finite word $w$ is {\em rich} if it contains exactly $|w| + 1$ distinct palindromes (or equivalently, if every prefix of $w$ has {\em Property $Ju$}). For example, $abac$ is rich since it is of length $4$ and contains the following five palindromes: $\empt$, $a$, $b$, $c$, $aba$. Naturally, an infinite word is rich if all of its factors are rich. For example, the periodic infinite words $a^\omega = aaa\cdots$ and $(ab)^\omega = ababab\cdots$ are clearly rich, whereas $(abc)^\omega = abcabacabc\cdots$ is not rich since it contains the non-rich word $abca$.

Droubay {\em et al.}~\cite{xDjJgP01epis} showed that all finite and infinite episturmian words are rich.  Specifically, they proved that if an infinite word has property~$Pi$ (and hence is epistandard -- see Theorem \ref{T:main}), then all of its prefixes have property~$Ju$. Consequently, any factor $u$ of an epistandard word (and hence, of an episturmian word) contains exactly $|u| + 1$ distinct palindromes, and is therefore rich (see Corollary~2 in \cite{xDjJgP01epis}).

Another special class of rich words the encompasses the episturmian words consists of Fischler's sequences with ``abundant palindromic prefixes''. These words were introduced and studied in~\cite{sF06pali, sF06pali2} in the context of Diophantine approximation. See also papers by Adamczewski and Bugeaud \cite{bAyB07pali, bAyB08tran} concerning the transcendence of certain real numbers whose sequences of partial quotients contain arbitrarily long palindromes. 

The theory of rich words has recently been further developed in a series of papers~\cite{aGjJ07rich, mBaDaGlZ08acon, aDaGlZ08rich, mBaDaGlZ08anew}.  In independent work, Ambro{\v{z}}, Frougny, Mas{\'a}kov{\'a}, and Pelantov{\'a}~\cite{pAzMePcF06pali} have considered the same class of words which they call {\em full words}, following the earlier work of Brlek, Hamel, Nivat, and Reutenauer in \cite{sBsHmNcR04onth}.

\subsection{Fractional powers \& critical exponent}

The study of {\em fractional powers} occurring in Sturmian words has been a topic of growing interest in recent times. See for instance \cite{jB99onth, vBcHlZ06init, dDdL02thei,  jJgP01frac, fMgP92repe, dV00stur}, as well as \cite{jJgP02epis, rRlZ00agen, aG05powe} for similar results concerning episturmian words and Arnoux-Rauzy sequences.

The following theorem extends the results in \cite{jJgP01frac} on fractional powers in Sturmian words. Throughout this section, we let $\bs$ denote an epistandard word with directive word $\Delta = x_1x_2\cdots \in \cA^\omega$ (as usual), and for all $n\geq 1$, we denote by $u_{n+1}$ the palindromic prefix $Pal(x_1\cdots x_n)$ of $\bs$ given in Theorem~\ref{T:main}. As in \cite{jJgP01frac}, we denote by $L(m)$ the length of the longest factor $v \in F(\bs)$ having period $m \in \NN$, and write $L(m) = em + r$, $e \in \NN^+$, $0 \leq r < m$. Given a finite or infinite word $w$, we denote by $w(i)$ (resp.~$w(i,j)$) the letter in position $i$ of $w$ (resp.~the factor of $w$ beginning at position $i$ and ending at position $j$).  

When $L(m) \geq 2m$,  all factors of $\bs$ having period $m$ and length $L(m)$ are equal to a palindrome $v$, and for $0 \leq i < e$, the word $v_i := v(1,im+r)$ is a palindromic prefix of $\bs$ by Lemma~4.1 in \cite{jJgP02epis}. Moreover, with the preceding notation, we have:

\begin{theo} {\em \cite[Theorem~4.2]{jJgP02epis}} 
Let $m$, $n \in \NN$ be such that $|u_n| < m \leq |u_{n+1}|$ and $\bs(1,m) = w$ is primitive with $\bs(m) = x$ occurring in $\bs(1,m-1)$. Then the following properties hold.
\begin{enumerate}
\item[i)] $L(m) \geq 2m$ if and only if $w = \mu_{x_1\cdots x_n}(x)$ and $x \in \Alphit(x_{n+1}x_{n+2}\cdots)$.
\item[ii)] Suppose $L(m) \geq 2$ and define $p = \max \{i \leq n ~|~ x_i = x\}$ and $t = \min\{j \in \NN^+ ~|~ x_{n+j} \ne x\}$. Then $u_{n+1} = w^tu_p$ is the longest prefix of $\bs$ having period $m$. Moreover, if $x \in \Alphit(x_{n+t+1}x_{n+t+2}\cdots)$, then $e = t + 1$; that is, $v = w^{t+1}u_p$, otherwise $e = t$ and $v = w^tu_p$. \qed
\end{enumerate} 
\end{theo}

\begin{rem} Let us mention a few noteworthy facts.
\begin{itemize}
\item Exponents of powers in $\bs$ are bounded if and only if exponents of letters in $\Delta$ are bounded~\cite{rRlZ00agen, jJgP02epis}.
\item Any Sturmian word has square prefixes and so do epistandard words \cite{jAjDmQlZ01tran, rRlZ00agen}. 
\item Any episturmian word has infinitely many prefixes of the form $uv^2$ with $|u|/|v|$ bounded above. 
\end{itemize}
The latter fact is readily deduced from the following result of Risley and Zamboni \cite{rRlZ00agen}.
\end{rem}

\begin{theo} {\em \cite[Prop.~I.3]{rRlZ00agen}} 
If $\bt$ is an Arnoux-Rauzy sequence, then there exists a positive number $\epsilon$ such that $\bt$ begins with infinitely many blocks of the form $UVVV'$, where $V'$ is a prefix of 
$V$ and $\min\{|V'|/|V|, |V|/|U|\} > \epsilon$. \qed
\end{theo}
\begin{note}
Such a result is motivated by transcendence issues; see for instance \cite{sFcM97tran}.
\end{note}

When $\bs$ is purely morphic, it is possible to give a rather explicit formula for the {\em critical exponent}: $\gamma = \limsup_{n\rightarrow\infty} L(m)/m$, as follows.

\begin{notation} Let $P$ be the function defined by $P(n) = \sup\{i < n ~|~ x_i = x_n\}$ if this integer exists, undefined otherwise. That is, if $x_n = a$, then $P(n)$ is the position of the right-most occurrence of the letter $a$ in the prefix $x_1x_2\cdots x_{n-1}$ of the directive word $\Delta = x_1x_2x_3\cdots \in \cAw$. 
\end{notation}

\begin{theo} {\em  \cite[Theorem~5.2]{jJgP02epis}}
Let $\bs$ be an $\cA$-strict epistandard word generated by a morphism with directive word $\Delta$ having period $q$.  Further, let $l \in \NN$ be maximal such that $y^l \in F(\Delta)$ for some letter $y$, and define $L = \{r, 0 \leq r < q ~|~ x_{r+1}= x_{r+2} = \cdots = x_{r+l}\}$ and $d(r) = r + q+1- P(r+q-1)$ for $0 \leq r < q$. Then the critical exponent for $\bs$ is given by
\[
  \gamma = l + 2 + \sup_{r\in L}\left\{\lim_{i\rightarrow\infty} |u_{r+iq+1-d(r)}|/|h_{r+iq}|\right\}.
\] 
Moreover, for any letter $x$ in $\bs$ the limit above can be obtained as a rational function with rational coefficients of the frequency $\alpha_x$ of this letter. \qed
\end{theo}
See also \cite{fMgP92repe, bTzW06some, dV00stur} for results on the critical exponent for the Fibonacci word, Tribonacci word, and Sturmian words, respectively.

\begin{exa} For the ever-so popular Fibonacci word $\bbf$, directed by $(ab)^\omega$, we have $q = 2$, $l=1$, $d(0)=d(1) =2$. Hence, since $|u_{n-1}|/|h_n|$ has limit $1/\varphi$ where $\varphi = (1+\sqrt{5})/2$ is the golden ratio, we obtain the well-known value $2 + \varphi$ for the critical exponent, originally proved by Mignosi and Pirillo~\cite{fMgP92repe}.

More generally, the $k$-bonacci word, directed by $(a_1a_2\cdots a_k)^\omega$, has critical exponent $2 + 1/(\varphi_k -1)$, where the {\em $k$-bonacci constant} $\varphi_k$ is the (unique) positive real root of the $k$-th degree monic polynomial 
$x^k - x^{k-1} - \cdots - x - 1$. 
%
\end{exa}

 \subsection{Frequencies} \label{SS:frequencies}

Let $w$ be a non-empty finite word. For any $v \in F(w)$, the {\em frequency} of $v$ in $w$ is $|w|_v/|w|$ where $|w|_v$ denotes the number of distinct occurrences of $v$ in $w$. The notion of frequency can be extended to infinite words in two ways, as follows.

\begin{defi} Suppose $v$ is a non-empty factor of an infinite word $\bx$. Then:
\begin{enumerate}
\item[i)] the frequency of $v$ in $\bx$ in the {\em weak sense} is $\lim_{n\rightarrow\infty} |w(1,n)|_v/n$ if this limit exists;
\item[ii)] $v$ has frequency $\alpha_v$ in $\bx$ in the {\em strong sense} if, for any sequence $(w_n)_{n\geq0}$ of factors of $\bx$ with increasing lengths, we have $\alpha_v = \lim_{n\rightarrow\infty} |w_n|_v/|w_n|$.
\end{enumerate}
\end{defi}

In a purely combinatorial way, Justin and Pirillo \cite[Section~6]{jJgP02epis} proved that any factor  occurring in an episturmian word has frequency in the strong sense. 

Wozny and Zamboni \cite{nWlZ01freq} also studied frequencies (in the weak sense) for Arnoux-Rauzy sequences. Using a reformulation of a vectorial division algorithm, originally introduced in \cite{rRlZ00agen}, they computed each allowable frequency of factors of the same length, as well as the number of factors with a given frequency. In particular, the authors of \cite{nWlZ01freq} gave simultaneous rational approximations of the frequencies by unreduced fractions having a common denominator. From this work, one recovers the results of Berth\'e  \cite{vB96freq} for Sturmian words in terms of {\em Farey} approximations arising from the continued fraction expansions of the frequencies of the letters. For instance, the frequencies of factors of the same length in a Sturmian word assume at most three values, which were explicitly given by Berth\'e \cite{vB96freq}, who also discovered that this result is in strong connection with the {\em three distance theorem} in Diophantine analysis.

\subsection{Return words} \label{S:return-words}

Let us now recall the notion of a {\em return word}, which was introduced independently by Durand \cite{fD98acha}, and Holton and Zamboni \cite{cHlZ99desc} when studying primitive substitutive sequences. 

\begin{defi}
Let $v$ be a recurrent factor of $\by \in \cA^\omega$, starting at positions $n_1 <n_2< n_3 \cdots$. Then each word $r_i = y_{n_i}y_{n_i+1}\cdots y_{n_{i+1} -1}$ is called a {\em return} to $v$ in $\by$. Moreover, $\by$ can be factorized in a unique way as $\by = y_1\cdots y_{n_1-1}r_1r_2r_3 \cdots$ where $r_1r_2r_3\cdots$, viewed as a word  on the $r_i$, is called the {\em derived word of $\by$ with respect to $v$}.
\end{defi}

That is, a return to $v$ in $\by$ is a non-empty factor of $\by$ beginning at an occurrence of $v$ and ending exactly before the next occurrence of $v$ in $\by$. Thus, if $r$ is a return to $v$ in $\by$, then $rv$ is a factor of $\by$ that contains exactly two occurrences of $v$, one as a prefix and one as a suffix. We call $rv$ a {\em complete return} to $v$ \cite{jJlV00retu}.  

Return words play an important role in the study of {\em minimal subshifts} in symbolic dynamics; see for instance \cite{fD98acha, Du2, Du1, FMN, Si}. In the context of episturmian words, such words have recently proven to be a useful tool in the study of quasiperiodicity (see Section~\ref{S:conclusion} for further details). This latest work made use of the following result of Justin and Vuillon \cite{jJlV00retu} which completely describes the returns to any factor of an epistandard word. In fact, their result actually characterizes return words in episturmian words (not just epistandard words) since, by uniform recurrence, the returns to any factor $v$ in an epistandard word $\bs$ are the same as the returns to $v$ as a factor of any episturmian word $\bt$ having the same set of factors as $\bs$.

\begin{theo} {\em \cite{jJlV00retu}} \label{T:return-words} 
Let $\bs$ be an epistandard word directed by $\Delta = x_1x_2x_3\cdots \in \cAw$ and consider any $v \in F(\bs)$. If $u_{n+1}$ is the shortest palindromic prefix of $\bs$ containing $v$ with $u_{n+1} = fvg$, then the returns to $v$ in $\bs$ are the words $f^{-1}\mu_{x_1\cdots x_n}(x)f$ where $x \in \Alphit(x_{n+1}x_{n+2}\cdots)$. Moreover, the corresponding complete returns to $v$ are the words $f^{-1}(u_{n+1}x)^{(+)}g^{-1}$ and the derived word of $\bs$ with respect to $v$ is given by $\bs^{(n)} = \mu_{x_1\cdots x_n}^{-1}(\bs)$. \qed
\end{theo}

\begin{note}It follows immediately that any factor of an $\cA$-strict episturmian word has exactly $|\cA|$ return words.
\end{note} 

Theorem~\ref{T:return-words} extends earlier work of Vuillon on return words in Sturmian words (see \cite{lV01acha}). In particular, Vuillon proved that Sturmian words are characterized by the property that any non-empty factor has exactly $2$ different return words in the given Sturmian word. However, contrary to what one might expect, such a property with $2$ replaced by a positive integer $k\geq 3$ does not characterize $k$-strict episturmian words. For instance, infinite words coding $3$-interval exchange transformations, which constitute a different generalization of Sturmian 
words to $3$-letter alphabets, are known to have the property that every factor has $3$ different return words  (see the work by Ferenczi, Holton, and Zamboni in \cite{sFcHlZ01stru}).

\section{Balance \& lexicographic order} \label{S:balance}


\subsection{$q$-Balance} \label{SS:q-balance}

\begin{defi} \label{D:q-balance}
A finite or infinite word is {\em $q$-balanced} if, for any two of its factors $u$, $v$ with $|u| = |v|$, we have
\[ 
||u|_x-|v|_x|\leq q  \quad \mbox{for any letter $x$},
\] 
i.e., the number of $x$'s in each of $u$ and $v$ differs by at most $q$.
\end{defi}

\begin{note} A $1$-balanced word is simply said to be {\em balanced}. \end{note}

The term `balanced' is relatively new; it appeared in \cite{jBpS93acha, jBpS94arem} (also see \cite[Chapter~2]{mL02alge}), and the notion itself dates back to \cite{gHmM40symb, eCgH73sequ}. In the pioneering work of Morse and Hedlund \cite{gHmM40symb}, balanced infinite words over a $2$-letter alphabet  were called `Sturmian trajectories' and belong to three classes: aperiodic Sturmian; periodic Sturmian; and infinite words that are ultimately periodic (but not periodic), called {\em skew words}. That is, the family of balanced infinite words consists of the (recurrent) Sturmian words and the (non-recurrent) skew infinite words, the factors of which are balanced. Skew words are ultimately periodic  suffixes of words of the form $\mu(a^pba^\omega)$, where $\mu$ is a {\em pure standard Sturmian morphism} and $p \in \NN$. For example, $aba^\omega$ and $\psi_b(aba^\omega) = bab(ba)^\omega$ are skew. See also \cite{rT96onco, rT98inte, aHrT00char, gP05mors} for further work on skew words.

\begin{remq} Nowadays, for most authors, only the aperiodic Sturmian words are considered to be `Sturmian'. However, from now on, we will use the term `Sturmian' to refer to both aperiodic and periodic Sturmian words. In the context of cutting sequences, the aperiodic (resp.~periodic) Sturmian words are precisely those with irrational slope (resp.~rational slope).
\end{remq}

It is important to note that a finite word is {\em finite Sturmian} (i.e., a factor of some Sturmian word) if and only if it is balanced \cite[Chapter 2]{mL02alge}. Accordingly, the balanced infinite words are precisely the infinite words whose factors are finite Sturmian. This concept was recently generalized in \cite{aGjJgP06char} by showing that the set of all infinite words whose factors are finite episturmian consists of the (recurrent) episturmian words and the (non-recurrent) {\em episkew} infinite words, as defined in the next section.

\subsection{Episkew words} \label{SS:episkew}

Inspired by the skew words of Morse and Hedlund \cite{gHmM40symb}, \emph{episkew words} were recently defined in \cite{aGjJgP06char} as non-recurrent infinite words, all of whose factors are (finite) episturmian. The following theorem gives a number of equivalent definitions of such words, similar to those for (recurrent) episturmian words. 

\begin{theo} {\em \cite{aGjJgP06char}} \label{T:skew} An infinite word $\bt$ with {\em Alph}$(\bt) = \cA$ is {\em episkew} if equivalently:
\begin{enumerate}
\item[i)] $\bt$ is non-recurrent and all of its factors are (finite) episturmian;
\item[ii)] there exists an infinite sequence $(\bt^{(i)})_{i\geq0}$ of non-recurrent infinite words and a directive word $x_1x_2x_3\cdots$ $(x_i \in \cA)$ such that $\bt^{(0)} = \bt$, $\ldots$~, $\bt'^{(i-1)} = \psi_{x_i}(\bt^{(i)})$, where  $ \bt'^{(i-1)}= \bt^{(i-1)}$ if $ \bt^{(i-1)}$ begins with $x_i$ and ${\bt^\prime}^{(i-1)} = x_i\bt^{(i-1)}$ otherwise;
\item[iii)]  there exists a letter $x \in \cA$ and an epistandard word $\bs$ on $\cA\setminus\{x\}$ such that $\bt = v \mu(\bs)$, where  $\mu$ is a pure epistandard morphism on $\cA$ and $v$ is a non-empty suffix of $\mu(\rev{\bs_p}x)$ for some $p \in \NN$.
\end{enumerate} \vspace{-0.2cm}
Moreover, $\bt$ is said to be {\em strict episkew} if $\bs$ is strict on $\cA\setminus\{x\}$, i.e., if each letter in $\cA\setminus\{x\}$ occurs infinitely often in the directive word $x_1x_2x_3\cdots$. \qed
\end{theo}
A simple example of an episkew word on more than two letters is the infinite word $c\bbf = cabaababa\cdots$ where $\bbf$ is the Fibonacci word and $c$ is a letter (see also Example~\ref{ex:episkew}).

Note that the episkew words on a $2$-letter alphabet are precisely the skew words. 
Certainly, in the Sturmian case, the word $\rev{\bs}_px\bs$ reduces to a word of the form $a^pba^\omega$.

\begin{rem}
Thanks to Richomme \cite{gR07priv}, episkew words actually have the following  simpler characterization: an infinite word $\bt$ is episkew if and only if $\bt = \varphi(x\bs)$ where $\bs$ is an  epistandard word, $x$ is a letter not occurring in $\bs$, and $\varphi$ is a pure episturmian morphism.   
\end{rem}

Episkew words were first alluded to (but not explicated) in the recent paper \cite{aG06acha}. 
Following that paper, these words showed up again in the study of inequalities characterizing finite and infinite episturmian words with respect to lexicographic orderings \cite{aGjJgP06char}. In fact, as detailed in the next section, episturmian words have similar extremal properties to Sturmian words.  See also \cite{jJgP02onac, oJlZ04char, gP03ineq,  gP05ineq, gP05mors, aG06acha, aGjJgP06char} for other work in this direction.

\subsection{Extremal properties}



Suppose the alphabet $\cA$ is totally ordered by the relation $<$. Then we 
can totally order $\cA^*$ by the \emph{lexicographic order} $\leq$ 
defined as follows. Given two words $u$, $v \in \cA^+$, we have $u
\leq v$ if and only if either $u$ is a prefix of $v$ or $u =
xau^\prime$ and $v = xbv^\prime$, for some $x$, $u^\prime$,
$v^\prime \in \cAstar$ and letters $a$, $b$ with $a < b$. This is
the usual alphabetic ordering in a dictionary. We write $u < v$ when $u\leq v$ and $u\ne v$, in which case we say that $u$ is (strictly) \emph{lexicographically smaller} than $v$. The notion of lexicographic order naturally extends to infinite words in $\cAw$. We denote by $\min(\cA)$ the smallest letter in $\cA$ with respect to the given lexicographic  order.

Let $w$ be a finite or infinite word over $\cA$ and let $k$ be a positive integer. We denote by $\min(w | k)$ (resp.~$\max(w | k)$) the lexicographically smallest (resp.~greatest) factor of $w$ of length $k$ for the given order (where $|w|\geq k$ if $w$ is finite). If $w$ is infinite, then it is clear that $\min(w | k)$ and $\max(w | k)$ are prefixes of the respective words $\min(w | k+1)$ and $\max(w | k+1)$. So we can define, by taking limits, the following two infinite words (see \cite{gP05ineq}):
\[
  \min(w) = \underset{k\rightarrow\infty} {\lim}\min(w | k) \quad \mbox{and} \quad 
  \max(w) = \underset{k\rightarrow\infty}{\lim}\max(w | k).
\]  
That is, to any infinite word $\bt$ we can associate two infinite words $\min(\bt)$ and $\max(\bt)$ such that any prefix of $\min(\bt)$ (resp.~$\max(\bt)$) is the lexicographically smallest (resp.~greatest)  amongst the factors of $\bt$ of the same length. 

For a finite word $w$ over $\cA$ and a given order on $\cA$, $\min(w)$ denotes $\min(w | k)$ where $k$ is maximal such that all $\min (w | j)$, $j= 1,2, \dots, k$, are prefixes of $\min (w | k)$. In the case $\cA= \{a,b\}$, $\max(w)$ is defined similarly (see \cite{aGjJgP06char}).

In 2003, Pirillo \cite{gP03ineq} (see also \cite{gP05ineq}) proved that, for infinite words $\bs$ on a $2$-letter alphabet $\{a,b\}$ with $a<b$, the inequality 
\begin{equation} \label{eq:characteristic}
a\bs \leq \min(\bs) \leq \max(\bs) \leq b\bs
\end{equation} 
characterizes {\em standard Sturmian words} (aperiodic and periodic). Actually, this result was known much earlier, dating back to the work of P.~Veerman \cite{pV86symb, pV87symb} in the mid 80's. Since that time, these `Sturmian inequalities' have been rediscovered numerous times under different guises, as discussed in the forthcoming survey paper \cite{jAaG06extr}.

Continuing his work in relation to inequality \eqref{eq:characteristic}, Pirillo \cite{gP05ineq} proved further that, in the case of an arbitrary finite alphabet $\cA$, an infinite word $\bs \in \cAw$ is epistandard if and only if, for any lexicographic order, we have 
\begin{equation} \label{eq:gP05ineq}
a\bs \leq \min(\bs) \quad \mbox{where $a = \min(\cA)$}. 
\end{equation} 
Moreover, $\bs$ is a strict epistandard word if and only if \eqref{eq:gP05ineq} holds with strict equality for any order~\cite{jJgP02onac}.

In a similar spirit, Glen, Justin, and Pirillo \cite{aGjJgP06char} recently established new characterizations of finite Sturmian and episturmian words via lexicographic orderings. As a consequence, they were able to characterize by lexicographic order all episturmian and episkew words. Similarly, they characterized by lexicographic order all balanced infinite words on a 2-letter alphabet; in other words, all Sturmian and skew infinite words, the factors of which are (finite) Sturmian. In the finite case:

 \begin{theo} {\em \cite{aGjJgP06char}} \label{T:aGjJgP06char}  A finite word $w$ on $\cA$ is episturmian if and only if there exists a finite word $u$ such that, for any lexicographic order,  
\begin{equation}au_{|m|-1} \le m \label{e2} \end{equation} where $m= \min(w)$ and $a = \min(\cA)$ for the considered order. \qed
\end{theo}

\newpage
\begin{exa} 
Consider the finite word $w= baabacababac$. For the different orders on $\{a,b,c\}$, we have
\begin{itemize}
\item $a<b<c$ or $a<c<b$: $\min(w) = aabacababac$,
\item $b<a<c$ or $b<c<a$: $\min(w) = babac$,
\item $c<a<b$ or $c<b<a$: $\min(w) = cababac$. 
\end{itemize}
It can be verified that a finite word $u$ satisfying \eqref{e2} must begin with $aba$ and one possibility is $u=abacaaaaaa$; thus $w$ is a finite episturmian word. 
\end{exa}
\begin{note} In the above example, any two orders with the same minimum letter give the same $\min(w)$, which is not true in general.
\end{note}

A corollary of Theorem \ref{T:aGjJgP06char} is the following new characterization of finite Sturmian words (i.e., finite balanced words).

\begin{cor} {\em \cite{aGjJgP06char}} \label{Cor:aGjJgP06char-finite}  A finite word $w$ on $\cA=\{a,b\}$, $a<b$, is not Sturmian (in other words, not  balanced)  if and only if there exists a finite word $u$ such that $aua$ is a prefix of $\min(w)$ and $bub$ is a prefix of $\max(w)$. 
\qed
\end{cor}

In the infinite case, the following characterization of all infinite words whose factors are finite episturmian  follows almost immediately from Theorem \ref{T:aGjJgP06char}. 

\begin{cor} {\em \cite{aGjJgP06char}} \label{Cor:aGjJgP06char-infinite}  An infinite word $\bt$ on $\cA$ is episturmian or episkew if and only if there exists an infinite word $\bu$ such that, for any lexicographic order, 
\begin{equation*}a\bu \le \min(\bt) \quad \mbox{where $a = \min(\cA)$.}  \vspace{-0.4cm} 
\end{equation*} \qed
\end{cor}

Consequently, an infinite word $\bs$ on $\{a,b\}$ ($a< b$) is balanced (i.e., Sturmian or skew) if and only if there exists an infinite word $\bu$ such that 
\begin{equation*}a\bu \le \min(\bs) \le \max(\bs) \le b\bu. 
\end{equation*}

Corollary~\ref{Cor:aGjJgP06char-infinite}  was recently refined in \cite{aG07orde} where it was shown  that, for any aperiodic episturmian word $\bt$, the infinite word $\bu$ (as given in the corollary) is the unique epistandard word with the same set of factors as $\bt$. As an easy consequence, we obtain the following  characterization  of strict episturmian words that are {\em infinite Lyndon words} (Theorem~\ref{T:Lyndon}).  Recall that a  non-empty finite word $w$ over $\cA$ is a {\em Lyndon word} if it is lexicographically smaller than all of its proper suffixes for the given order $<$ on $\cA$. Equivalently, $w$ is the lexicographically smallest primitive word in its conjugacy class; that is, $w < vu$ for all non-empty words $u$, $v$ such that $w = uv$. The first of these definitions extends to infinite words: an infinite word over $\cA$ is an {\em infinite Lyndon word} if and only if it is (strictly) lexicographically smaller than all of its proper suffixes for the given order on $\cA$. That is, a finite or infinite word $w$ is a Lyndon word if and only if $w < \TT^i(w)$ for all $i > 0$. 

Assuming that $|\cA| > 1$ (since there are no Lyndon words on a $1$-letter alphabet), we have:

\begin{theo} {\em \cite{aG07orde}} \label{T:Lyndon}
An $\cA$-strict episturmian word $\bt$ is an infinite Lyndon word if and only if $\bt = a\bs$ where $a = \min(\cA)$ for the given order on $\cA$ and $\bs$ is an (aperiodic) $\cA$-strict epistandard word. Moreover, if $\Delta = x_1x_2\cdots \in \cAw$ is the directive word of $\bs$, then $\bt = a\bs$ is the unique episturmian word in the subshift of $\bs$  directed by the spinned version of $\Delta$ having all spins $1$, except when $x_i = a$.  \qed
\end{theo}
The above theorem is actually a generalization of a result on (aperiodic) Sturmian words given by Borel and Laubie \cite{jBfL93quel} (see also \cite{gR07conj}). 

Let $\cA = \{a_1, \ldots, a_m\}$ be an alphabet ordered by $a_1 < a_2 < \cdots < a_m$. Then Theorem~\ref{T:Lyndon} says that an $\cA$-strict episturmian word $\bt$ is an infinite Lyndon word if and only if the (normalized) directive word of $\bt$ belongs to $\{a_1,\bar a_2, \ldots, \bar a_m\}^\omega$. This can be reformulated as a generalization of Proposition~6.4 in \cite{fLgR07quas}:

\begin{cor} {\em \cite{aG07orde}} \label{Cor:decomposable2} An $\cA$-strict episturmian word $\bt$ is an infinite Lyndon word if and only if it can be infinitely decomposed over the set of morphisms $\{\psi_a, \bar\psi_x \mid x \in \cA\setminus\{a\}\}$ where $a = \min(\cA)$ for the given order on~$\cA$. \qed
\end{cor} 

We observe that, contrary to the fact that there exists $|\cA|!$ possible orders of a finite alphabet $\cA$, Theorem~\ref{T:Lyndon} shows that there exist exactly $|\cA|$ infinite Lyndon words in the subshift of a given $\cA$-strict epistandard word $\bs$ (when $|\cA| >1$). That is, for any order with $\min(\cA) = a$, the subshift of $\bs$ contains a unique infinite Lyndon word beginning with $a$, namely $a\bs$.  

\begin{exa}
With $\Delta = (abc)^\omega$, the spinned versions $(a\bar b \bar c)^\omega$, $(\bar a b \bar c)^\omega$, $(\bar a \bar b c)^\omega$ and their `opposites' (obtained by exchange of spins): $(\bar a b c)^\omega$, $(a\bar b c)^\omega$, $(ab\bar c)^\omega$ direct episturmian words in the subshift of the Tribonacci word $\br$.  Only the first three of these spinned infinite words direct episturmian Lyndon words: $a\br$, $b\br$, $c\br$, respectively.
\end{exa}

The above results on strict episturmian Lyndon words have very recently been generalized to all episturmian words by Glen, Lev\'e, and Richomme~\cite{aGfLgR07quas}, as follows.

\begin{theo} \label{T:Lyndon-episturmian} {\em \cite{aGfLgR07quas}} 
Let $\cA = \{a_1, \ldots, a_m\}$ be an alphabet ordered by $a_1 < a_2 < \cdots < a_m$ and, for $1 \leq i \leq m$, let $\cB_i = \{a_i, \ldots, a_m\}$.  An episturmian word $\bt$ is an infinite Lyndon word if and only if there exists an integer $j$ such that $1 \leq j <m$ and the (normalized) directive word of $\bw$ belongs to:
$$(\bar\cB_2^*a_1)^*\cdots 
(\bar\cB_j^*a_{j-1})^*(\bar\cB_{j+1}^*a_j)^*(\bar\cB_{j+1}^+\{a_j\}^+)^\omega.$$ \qed
\end{theo}
\begin{note} In the above theorem, the word \emph{normalized} appears between brackets since one can easily verify from Theorem~\ref{T:uniqueDirective} that a spinned infinite word of the given form is the unique directive word of exactly one episturmian word. 
\end{note}

\begin{exa} \cite{aGfLgR07quas} 
Let $\cA=\{a,b,c,d\}$. Then the spinned infinite word $(\bar b \bar c a)(\bar d \bar c b)^2 (\bar d c c)^\omega$ directs a Lyndon episturmian word, and so does $aa(\bar d c)^\omega$, but $\bar c a \bar b a \bar d c d^\omega$ does not (since this spinned word directs a periodic word).
\end{exa}

\begin{rem} \label{R:uniqueDirective-Lyndon}
Theorems \ref{T:uniqueDirective} and \ref{T:Lyndon-episturmian} show that any episturmian Lyndon word has a unique spinned directive word, but the converse is not true. For example, the regular wavy word $(a\bar b c)^\omega$ is the unique directive word of the strict episturmian word: 
\[
\lim_{n\rightarrow\infty}{\mu_{a\bar b c}^n(a)} = acabaabacabacabaabaca\cdots 
\]
which is clearly not an infinite Lyndon word by Theorem~\ref{T:Lyndon-episturmian} and also by the fact that $acabaaw$ is not a Lyndon word for any order on $\{a,b,c\}$ and for any word $w$. 
\end{rem}

A key tool used in the proof of Theorem~\ref{T:Lyndon-episturmian} was the following result of Richomme, which characterizes episturmian morphisms that preserve Lyndon words. A morphism $f$ is said to preserve finite (resp.~infinite) Lyndon words if for each finite (resp.~infinite) Lyndon word $w$, $f(w)$ is a finite (resp.~infinite) Lyndon word. 

\begin{theo} {\em \cite{gR03lynd, gR04onmo}} \label{P:Lyndon-morphisms}
Let $\cA = \{a_1, \ldots, a_m\}$ be an alphabet ordered by $a_1 < a_2 < \cdots < a_m$. Then the following assertions are equivalent for an episturmian morphism:
\begin{itemize}
\item $f$ preserves finite Lyndon words;
\item $f$ preserves infinite Lyndon words;
\item $f \in ({\bar\Psi}_{\{a_2, \ldots, a_m\}}^*\psi_{a_1})^*\{\bar\Psi_{a_m}\}^*$ where $\bar\Psi_{\cA} = \{\bar\psi_x ~|~ x \in \cA\}$. \qed
\end{itemize}
\end{theo}

\subsection{Imbalance}

We now return our attention to the notion of balance. 

Episturmian words on three or more letters are generally unbalanced in the sense of $1$-balance, except, of course, for those on a $2$-letter alphabet, which are precisely the (periodic and aperiodic)  Sturmian words. In fact, Cassaigne, Ferenczi, and Zamboni \cite{jCsFlZ00imba} have proved, by construction, that there exists an episturmian word that is not $q$-balanced for any $q$.  Note, however, that the Tribonacci word is $2$-balanced, for example. More generally, it can be shown by induction that the {\em $k$-bonacci word}, directed by $(a_1a_2\cdots a_k)^\omega$, is $(k-1)$-balanced. Even  further, one can prove that any {\em linearly recurrent} strict episturmian word (or Arnoux-Rauzy sequence) is $q$-balanced for some $q$. Linearly recurrent Arnoux-Rauzy sequences were completely described in \cite{rRlZ00agen, jCnC06fonc}; they are the strict episturmian words for which each letter $x$ occurs in $\Delta$ with bounded gaps. 

Using their main result on return words (Theorem \ref{T:return-words}), Justin and Vuillon \cite{jJlV00retu} proved that episturmian words do in fact satisfy a kind of balance property. Specifically:

\begin{theo} {\em \cite[Theorem~5.2]{jJlV00retu}} Let $\bs \in \cAw$ be an epistandard word and let $\{d,e\}$ be a $2$-letter subset of $\cA$. Then, for any $u$, $v \in F(\bs) \cap \{d,e\}^*$ with $|u| = |v|$, we have $||u|_d - |v|_d|\leq 1$. \qed
\end{theo}
This property of episturmian words reduces to the balance property of Sturmian words when $\cA$ is a $2$-letter alphabet (in which case it is characteristic); however, the property is far from being characteristic when $\cA$ consists of more than two letters.

More recently, Richomme \cite{gR07aloc} also proved that episturmian words and Arnoux-Rauzy sequences can be characterized via a nice `local balance property'.  
That is:

\begin{theo} \label{T:gR07aloc} {\em \cite{gR07aloc}} For a recurrent infinite word $\bt \in \cAw$, the following assertions are equivalent:
\begin{enumerate}
\item[i)] $\bt$ is episturmian; 
\item[ii)] for each factor $u$ of $\bt$, there exists a letter $a$ such that $\cA u\cA\cap F(\bt) \subseteq au\cA\cup\cA ua$;
\item[iii)] for each palindromic factor $u$ of $\bt$, there exists a letter $a$ such that $\cA u\cA\cap F(\bt) \subseteq  au\cA \cup \cA ua$. \qed
\end{enumerate}
\end{theo} 

Roughly speaking, the above theorem says that for any factor $u$ of a given episturmian word $\bt$, there exists a unique letter $a$ such that  every occurrence of $u$ in $\bt$ is immediately preceded or followed by $a$ in $\bt$. When $|\cA| = 2$, property $ii)$ of Theorem~\ref{T:gR07aloc} is equivalent to the definition of balance. Indeed, Coven and Hedlund \cite{eCgH73sequ} stated that an infinite word $\bw$ over $\{a, b\}$ is not balanced if and only if there exists a palindrome $u$ such that $aua$ and $bub$ are both factors of $\bw$. As pointed out in \cite{gR07aloc}, this property can be rephrased as follows: an infinite word $\bw$ is Sturmian if and only if $\bw$ is aperiodic and, for any factor $u$ of $\bw$, the set of factors belonging to $\cA u \cA$ is a subset of $au\cA \cup \cA u a$ or a subset of $bu\cA \cup \cA ub$.

\subsection{Fraenkel's conjecture}

As discussed previously, the {\em recurrent} balanced infinite words on two letters are exactly the Sturmian words (aperiodic and periodic). A natural question to ask is then: ``What are the balanced recurrent infinite words on more than two letters?'' In this direction, Paquin and Vuillon \cite{gPlV07acha}  recently characterized the balanced episturmian words  by classifying these words into three families, as follows.

\begin{theo} {\em \cite{gPlV07acha}}  \label{T:gPlV07acha} 
Any balanced standard episturmian sequence $\bs$ on a $k$-letter alphabet $\cA_k =\{1,2,\ldots, k\}$, $k\geq 3$,  belongs to one of the following three families (up to letter permutation):  
\begin{enumerate}
\item[i)] $\bs = p(k - 1)p(kp(k - 1)p)^\omega$, with $p = Pal(1^n 2 \cdots (k - 2))$; 
\item[ii)] $\bs =  p(k - 1)p (kp(k - 1)p)^\omega$, with $$p = Pal(123 \cdots (k - \ell - 1)1(k - \ell) \cdots (k - 2));$$
\item[iii)] $\bs = [Pal(123\cdots k)]^\omega$. \qed
\end{enumerate} 
\end{theo}

The importance of the above result lies in the fact that it supports {\em Fraenkel's conjecture} \cite{aF73comp}:  a problem that arose in a number-theoretic context and has remained unsolved for over thirty years. Fraenkel conjectured that, for a fixed $k \geq 3$, there is only one covering of $\ZZ$ by $k$ {\em Beatty sequences} of the form $(\lfloor \alpha n + \beta \rfloor)_{n\geq1}$, where $\alpha$, $\beta$ are real numbers. A combinatorial interpretation of this conjecture may be stated as follows (taken from  \cite{gPlV07acha}). Over a $k$-letter alphabet with $k\geq 3$, there is only one recurrent balanced infinite word, up to letter permutation  and shifts, that has mutually distinct letter frequencies. This supposedly unique infinite word is called {\em Fraenkel's sequence} and is given by $(F_k)^\omega$ where the {\em Fraenkel words} $(F_i)_{i\geq 1}$ are defined recursively by $F_1=1$ and  $F_{i}=F_{i-1}iF_{i-1}$ for  all $i\geq 2.$  (Note that $F_k = Pal(12\cdots k)$.) For further details, see for instance \cite{gPlV07acha, rT00frae} and references therein. 

Amongst the classes of balanced episturmian words given in Theorem~\ref{T:gPlV07acha}, only one class has mutually distinct letter frequencies and, up to letter permutation and shifts, corresponds to Fraenkel's sequence. That is:

\begin{theo}[Paquin-Vuillon \cite{gPlV07acha}] Suppose $\bt$ is a balanced episturmian word with $\Alphit(\bt) = \{1,2,\ldots, k\}$, $k \geq 3$. If $\bt$ has mutually distinct letter frequencies, then up to letter permutation, $\bt$ is a shift of~$(F_k)^\omega$. \qed
\end{theo}

More recently, it was proved in \cite{aGjJ07rich} that any recurrent balanced {\em rich} infinite word is necessarily episturmian, and hence such words obey Fraenkel's conjecture (recall that rich words were defined Section~\ref{SS:rich}).

\begin{remq}
An interesting known fact (e.g., see \cite{pH00suit}) is that any balanced recurrent infinite word $\bx$ on $k\geq 3$ letters having mutually distinct letter frequencies is necessarily periodic.  Certainly, the image of $\bx$ under any morphism of the form: ($a \mapsto a$, \mbox{other}~$x \mapsto b$) is a Sturmian word. If, for one letter, the corresponding  
Sturmian word is aperiodic (i.e., $\bx$ has irrational slope as a cutting sequence), then we meet impossibility;  thus rather easily $\bx$ must be periodic.
\end{remq}

\section{Concluding remarks} \label{S:conclusion}

In closing, we mention a number of very recent works involving episturmian words.

\begin{description}
\item[Rigidity:]
Krieger \cite{dK07onst} has shown that any strict purely morphic epistandard word $\bs$ is {\em rigid}. That is, all of the morphisms that generate $\bs$ are powers of the same unique (epistandard) morphism. Krieger also showed that a certain class of `ultimately strict' purely morphic epistandard words are not rigid, but it remains an open question as to whether or not all strict morphic episturmian words are rigid.

\item[Quasiperiodicity:]
A finite or infinite word $w$ is said to be {\em quasiperiodic} if there exists a word $u$ (with $u \ne w$ for finite $w$) such that the occurrences of $u$ in $w$ entirely cover $w$, i.e., every position of $w$ falls within some occurrence of $u$ in $w$. Such a word $u$ is called a {\em quasiperiod} of $w$. For example, the word $w=abaababaabaababaaba$ has quasiperiods $aba$, $abaaba$, $abaababaaba$. 

In the last fifteen years, quasiperiodicity and coverings of finite words has been extensively studied (see \cite{aAmC01stri} for a brief survey on quasiperiodicity in `strings').  Quasiperiodic finite words were first introduced by Apostolico and Ehrenfeucht in \cite{aAaE93effi}. The notion was later extended to infinite words by Marcus \cite{sM04quas} who opened some questions, particularly concerning quasiperiodicity of Sturmian words. After a brief answer to some of these questions in \cite{fLgR04quas}, the Sturmian case was fully studied by Lev\'e and Richomme~\cite{fLgR07quas} who proved that a Sturmian word is non-quasiperiodic if and only if it is an infinite Lyndon word. The study of quasiperiodicity in Sturmian words was very recently extended to episturmian words by Glen, Lev\'e, and Richomme~\cite{aG07orde, aGfLgR07quas, fLgR07quasB}, who have completely described all of the quasiperiods of an episturmian word, yielding a characterization of quasiperiodic episturmian words in terms of their directive words. They have also characterized episturmian morphisms that map any word onto a quasiperiodic one. These results show that, unlike the Sturmian case, there exist non-quasiperiodic episturmian words that are not infinite Lyndon words. Key tools used in the study of quasiperiodicity in episturmian words were episturmian morphisms, normalized directive words (recall Theorem~\ref{T:normalisation}), and the following equivalent definition of quasiperiodicity in terms of return words introduced by Glen in \cite{aG07orde}: a finite word $v$ is a quasiperiod of an infinite word $\bw$ if and only if $v$ is a recurrent prefix of $\bw$ such that all of the returns to $v$ in $\bw$ have length at most $|v|$. 

In \cite{tM05illu}, Monteil proved that any Sturmian subshift contains a {\em multi-scale quasiperiodic word}, i.e., an infinite word having infinitely many quasiperiods. A shorter proof of this fact was provided in \cite{fLgR07quas} and this result has also been proven true for episturmian words in \cite{aGfLgR07quas}. 

For more recent work on quasiperiodicity, see for instance \cite{tM05illu, tMsM08quas}.

\item[$\theta$-episturmian words:] Recall that an infinite word is episturmian if and only if its set of factors is closed under reversal and it has at most one left special factor of each length. With this definition in mind, Bucci, de Luca, De Luca, and Zamboni \cite{mBaDaDlZ07onso, mBaDaDlZ07ondi} have recently introduced and studied a further extension of episturmian words in which the reversal operator is replaced by an arbitrary {\em  involutory antimorphism} (i.e., a map $\theta: \cA^* \rightarrow \cA^*$ such that $\theta^2 =$ Id and $\theta(uv) = \theta(v)\theta(u)$ for all $u$, $v \in \cA^*$).  More precisely, an infinite word over $\cA$ is said to be {$\theta$-episturmian} if  it has at most one left special factor of each length and its set of factors is closed under an involutory antimorphism $\theta$ of the free monoid $\cA^*$. Generalizing even further, {\em $\theta$-episturmian words with seed} are obtained by requiring the condition on special factors only for sufficiently large lengths (see \cite{mBaDaDlZ07ondi}). 
\end{description}

\bigskip
\noindent {\large\bf Acknowledgements.} The authors would like to thank Jean Berstel and Pierre Arnoux for their helpful comments on a preliminary version of this paper. Many thanks also to the two anonymous referees whose thoughtful suggestions helped to improve the paper.


\small

\begin{thebibliography}{99}



\bibitem{bA03bala} B.~Adamczewski, Balances for fixed points of primitive substitutions, {\it Theoret.  Comput. Sci.} 307 (2003) 47--75. 

\bibitem{bAyB07pali} B.~Adamczewski,Y. Bugeaud, Palindromic continued fractions, {\it Ann. Inst. Fourier (Grenoble)} 57 (2007) 1557--1574.

\bibitem{bAyB08tran} B.~Adamczewski, Y.~Bugeaud, Transcendence measure for continued fractions involving repetitive or symmetric patterns, {\em J. Eur. Math. Soc.}, to appear. 

\bibitem{pAvB98three} P.~Alessandri, V.~Berth\'e, Three distance theorems and combinatorics on words, {\it Enseign. Math.} 44 (1998) 103--132.

\bibitem{jAjDmQlZ01tran}
J.-P.~Allouche, J.L.~Davison, M.~Queff\'{e}lec, L.Q.~Zamboni, 
  Transcendence of {S}turmian or morphic continued fractions, {\it J. Number Theory} 91 (2001) 39--66.

\bibitem{jAaG06extr} J.-P.~Allouche, A.~Glen, Extremal properties of (epi)sturmian sequences and distribution modulo $1$, in preparation.

\bibitem{jAjS03auto} J.-P.~Allouche, J.~Shallit,  \emph{Automatic Sequences: Theory,
               Applications, Generalizations},  \emph{Cambridge University Press}, UK, 2003.
               
\bibitem{pAzMePcF06pali} P.~Ambro{\v{z}}, C.~Frougny, Z.~Mas{\'a}kov{\'a}, E.~Pelantov{\'a},  Palindromic complexity of infinite words associated with 
  simple {P}arry numbers, {\it Ann. Inst. Fourier (Grenoble)} 56 (2006) 2131--2160.
               
\bibitem{aAmC01stri} A.~Apostolico, M.~Crochemore, String pattern matching for a deluge survival kit, in: {\it Handbook of Massive Data Sets, Massive Comput.}, vol.~4, Kluwer Acad. Publ., Dordrecht, 2002. 
               
\bibitem{aAaE93effi} A.~Apostolico, A.~Ehrenfeucht, Efficient detection of quasiperiodicities in strings, {\it Theoret. Comput. Sci.} 119 (1993) 247--265.               
               
\bibitem{pAsI01piso} P.~Arnoux, S.~Ito, Pisot substitutions and Rauzy fractals, {\it Bull. Belg. Math. Soc. Simon Stevin} {8} (2001) 181--207.               

\bibitem{pAgR91repr} P.~Arnoux, G.~Rauzy, Repr\'{e}sentation 
g\'{e}om\'{e}trique de suites de complexit\'{e} $2n+1$, 
\emph{Bull. Soc. Math. France} {119} (1991) 199--215.

\bibitem{yB95comp} Yu.~Baryshnikov, Complexity of trajectories in rectangular billiards, {\it  Comm. Math. Phys.} 174 (1995) 43--56.


\bibitem{jB99onth} J.~Berstel, On the index of Sturmian words, in: 
\emph{Jewels Are Forever}, Springer-Verlag, Berlin, 1999, pp.~287--294.

\bibitem{jB02rece} J.~Berstel, Recent results on extensions of {S}turmian words, {\it Internat. J. Algebra Comput.} {12} (2002) 371--385. 

\bibitem{jBpS93acha} J.~Berstel, P.~S\'{e}\'{e}bold, A characterization of {Sturmian morphisms}, in:  {\em Borzyszkowski, A.M. and Sokolowski, S. $($Eds.$)$, Mathematical Foundations of Computer Science 1993, Lecture Notes in Computer Science}, vol.~711, {\em Springer-Verlag}, Berlin, 1993, pp. 281--290.

\bibitem{jBpS94arem} J.~Berstel, P.~S\'{e}\'{e}bold, A remark on morphic {S}turmian words,  {\it Theor. Inform. Appl.} {28} (1994) 255--263.

\bibitem{jBlV02codi} J.~Berstel, L.~Vuillon, Coding rotations on intervals, {\it Theoret. Comput. Sci.} {281} (2002) 99--107. 

\bibitem{vB96freq} V.~Berth\'e, Fr\'equences des facteurs des suites sturmiennes, {\em 
Theoret. Comput. Sci.} {165} (1996) 295--309. 

\bibitem{vB01auto} V.~Berth\'{e}, Autour du syst\`{e}me de num\'{e}ration {d'Ostrowski}, {\it Bull. Belg. Math. Soc. Simon Stevin} {8} (2001) 209--239.

\bibitem{vBhEsIhR07onsu}
V.~Berth\'{e}, H.~Ei, S.~Ito, and H.~Rao, On substitution invariant {S}turmian words: an application of {R}auzy fractals, {\it Theoret. Inform. Appl.} 41 (2007) 329--349.

\bibitem{vBsFlZ05inte} V.~Berth\'e, S.~Ferenczi, L.Q.~Zamboni, Interactions between dynamics, arithmetics and combinatorics: the good, the bad, and the ugly, in: {\em Algebraic and topological dynamics}, {\em Contemp. Math.}, vol.~385, Amer. Math. Soc., Providence, RI, 2005, pp.~333--364. 

\bibitem{vBcHlZ06init} V.~Berth\'{e}, C.~Holton, L.Q.~Zamboni, Initial 
powers of {S}turmian sequences, {\it Acta Arith.} {122} (2006) 315--347.

\bibitem{jBfL93quel} J.-P.~Borel, F.~Laubie, Quelques mots sur la droite 
projective r\'eelle, {\it J. Th\'eor. Nombres Bordeaux\,} 5 (1993) 23--51.

\bibitem{jBcR05pali} J.-P.~Borel, C.~Reutenauer, Palindromic factors of billiard words, 
  {\it Theoret. Comput. Sci.} {340} (2005) 334--348. 
  
\bibitem{sBsHmNcR04onth} S.~Brlek, S.~Hamel, M.~Nivat, C.~Reutenauer,  On the palindromic complexity of infinite words, {\it Internat. J. Found. Comput. Sci.} { 15} (2004) 293--306.

\bibitem{mBaDaDlZ07onso} M.~Bucci, A.~de~Luca, A.~De~Luca, L.Q.~Zamboni, On some problems related to palindrome closure, {\it Theor. Inform. Appl.} (in press), doi:10.1051/ita:2007064.

\bibitem{mBaDaDlZ07ondi} M.~Bucci, A.~de~Luca, A.~De~Luca, L.Q.~Zamboni, On different generalizations of episturmian words, {\it Theoret. Comput. Sci.} 393 (2008) 23--36.

\bibitem{mBaDaGlZ08acon} M.~Bucci, A.~De~Luca, A.~Glen, L.Q.~Zamboni, A connection between palindromic and factor complexity using return words, {\it Adv. in Appl. Math.} (in press), doi:10.1016/j.aam.2008.03.005.

\bibitem{mBaDaGlZ08anew} M.~Bucci, A.~De~Luca, A.~Glen, L.Q.~Zamboni, A new characteristic property of rich words, Preprint, 2008, arXiv:0807.2303.

\bibitem{jC98sequ} J.~Cassaigne, Sequences with grouped factors, in: {\it Developments in
  Language Theory III}, Aristotle University of Thessaloniki, 1998, pp.~211--222.

\bibitem{jCnC06fonc} J.~Cassaigne, N.~Chekhova, Fonctions de r\'ecurrence des
  suites d'{A}rnoux-{R}auzy et r\'eponse \`a une question de {M}orse et {H}edlund, {\it Ann. Inst. Fourier (Grenoble)} {56} (2006) 2249--2270.  

\bibitem{jCsFlZ00imba} J.~Cassaigne, S.~Ferenczi, L.Q.~Zamboni, Imbalances in
  {Arnoux-Rauzy} sequences, {\it Ann. Inst. Fourier (Grenoble)} {50} (2000) 1265--1276.
  
\bibitem{mCfMaR99fine} M.G.~Castelli, F.~Mignosi, A.~Restivo, Fine and {W}ilf's theorem for
  three periods and a generalization of {S}turmian words, {\it Theoret. Comput. Sci.} {218} (1999) 83--94. 
  
\bibitem{nCpHaM01prop} N.~Chekhova, P.~Hubert, A.~Messaoudi, Propri\'et\'es combinatoires, ergodiques et arithm\'etiques de la substitution de Tribonacci, {\it J. Th\'eor. Nombres Bordeaux\,} 13 (2001) 371--394.
  
\bibitem{eC74sequ} E.M.~Coven, Sequences with minimal block growth II, {\it Math. Systems Theory} 8 (1974) 376--382.  
  
\bibitem{eCgH73sequ} E.M.~Coven, G.A.~Hedlund, Sequences with minimal block growth, \emph{Math. Systems Theory} {7} (1973) 138--153.  

\bibitem{dCwMaPpS93subs} D.~Crisp, W.~Moran, A.~Pollington, P.~Shiue, Substitution invariant cutting sequences, {\it J. Th\'eor. Nombres Bordeaux} 5 (1993) 123--137.
  
\bibitem{dDdL02thei} D.~Damanik, D.~Lenz, The index of Sturmian sequences, 
\emph{European J. Combin.} {23} (2002) 23--29.  

\bibitem{dDlZ03comb} D.~Damanik, L.Q.~Zamboni, Combinatorial properties of
  {A}rnoux-{R}auzy subshifts and applications to {S}chr\"odinger operators, 
  {\it Rev. Math. Phys.} {15} (2003) 745--763. 
  
\bibitem{aD97stur} A.~de~Luca, Sturmian words: structure,
combinatorics and their arithmetics, \emph{Theoret. Comput. Sci.}
{183} (1997) 45--82. 

\bibitem{aDaGlZ08rich} A.~de~Luca, A.~Glen, L.Q.~Zamboni, Rich, Sturmian, and trapezoidal words, {\em Theoret. Comput. Sci.} (in press), doi:10.1016/j.tcs.2008.06.009. 
  
\bibitem{xDjJgP01epis} X.~Droubay, J.~Justin, G.~Pirillo,
Episturmian words and some constructions of de Luca and Rauzy,
\emph{Theoret. Comput. Sci.} {255} (2001) 539--553. 

\bibitem{xDgP99pali} X.~Droubay, G.~Pirillo, Palindromes and {S}turmian words, {\it Theoret. 
  Comput. Sci.} {223} (1999) 73--85. 


\bibitem{fD98acha} F.~Durand, A characterization of substitutive sequences using return words, {\em Discrete Math.} {179} (1998) 89--101. 

  
  \bibitem{Du2} F.~Durand, A generalization of Cobham's theorem, {\em Theory Comput. Syst.} {31} (1998) 169--185.
  
  \bibitem{Du1} F.~Durand, Linearly recurrent subshifts have a finite number of non-periodic subshift factors, {\em Ergodic Theory Dynam. Systems} {19} (1999) 953--993.
  


\bibitem{iF06alit} I.~Fagnot, A little more about morphic Sturmian words, {\it Theor. Inform. Appl.} {40}  (2006) 511--518. 

\bibitem{iFlV02gene} I.~Fagnot, L.~Vuillon, Generalized balances in {S}turmian words, {\it Discrete Appl.  Math.} 121 (2002) 83--101.


\bibitem{sF99comp} S.~Ferenczi, Complexity of sequences and dynamical systems, {\em Discrete Math.} {206} (1999) 145--154. 

\bibitem{sFcHlZ01stru} S.~Ferenczi, C.~Holton, L.Q.~Zamboni, Structure of three interval exchange transformations {I}. {A}n arithmetic study, {\it Ann. Inst. Fourier (Grenoble)} {51} (2001) 861--901.

\bibitem{sFcM97tran} S.~Ferenczi, C.~Mauduit, Transcendence of numbers with a low
  complexity expansion, {\it J. Number Theory} 67 (1997) 146--161.


  \bibitem{FMN} S.~Ferenczi, C.~Mauduit, A.~Nogueira, Substitutional dynamical systems: algebraic characterization of eigenvalues, {\em Ann. Sci. \'Ecole Norm. Sup.} {29} (1995) 519--533.


\bibitem{sF06pali} S.~Fischler, Palindromic prefixes and episturmian words,  
  {\it J. Combin. Theory Ser. A} {113} (2006) 1281--1304. 
  
  
\bibitem{sF06pali2}
S.~Fischler, Palindromic prefixes and diophantine approximation, {\it Monatsh. Math.} 151 (2007) 11--37.
  
\bibitem{aF73comp} A.S.~Fraenkel, Complementing and exactly covering sequences. {\it J. Combinatorial Theory Ser. A} {14} (1973) 8--20. 

\bibitem{aG06onst} A.~Glen, On {Sturmian} and episturmian words, and related topics, Ph.D. Thesis, The University of Adelaide, Australia, April 2006. 


\bibitem{aG07orde} A.~Glen. Order and quasiperiodicity in episturmian words. in: {\em Proceedings of the $6$th International Conference on Words}, Marseille, France, September 17-21, 2007, pp.~144--158.

\bibitem{aG05powe} A.~Glen, Powers in a class of $\mathcal{A}$-strict standard 
               episturmian words,  {\em Theoret. Comput. Sci.} {380} (2007) 330--354. 
               
\bibitem{aG06acha}
A.~Glen, A characterization of fine words over a finite alphabet,
  {\em Theoret. Comput. Sci.} 391 (2008) 51--60. 


\bibitem{aGjJ07rich} A.~Glen, J.~Justin, S.~Widmer, L.Q.~Zamboni, Palindromic richness, {\it European J. Combin.} (in press), doi:10.1016/j.ejc.2008.04.006.

\bibitem{aGjJgP06char} A.~Glen, J.~Justin, G.~Pirillo, Characterizations 
of finite and infinite episturmian words via lexicographic orderings, \emph{European J. Combin.} 29 (2008) 45--58. 

\bibitem{aGfLgR08dire} A.~Glen, F.~Lev\'e, G.~Richomme, Directive words of episturmian words: equivalences and normalization, Preprint, 2008, arXiv:0802.3888.

\bibitem{aGfLgR07quas} A.~Glen, F.~Lev\'e, G.~Richomme, Quasiperiodic and Lyndon episturmian words, {\it Theoret. Comput. Sci.}, to appear.

\bibitem{eG07repr} E.~Godelle, Repr\'esentation par des transvections des groupes dÕartin-tits, {\em Group, Geometry and Dynamics} 1 (2007) 111--133. 


\bibitem{aHrT00char} A.~Heinis, R.~Tijdeman, Characterisation of asymptotically {S}turmian sequences, \emph{Publ. Math. Debrecen} {56} (2000) 415--430. 

\bibitem{cHlZ99desc} C.~Holton, L.Q.~Zamboni, Descendants of primitive substitutions, {\em Theory Comput. Syst.} {32} (1999) 133--157. 

\bibitem{pH00suit} P. Hubert, Suites \'equilibr\'es, {\it Theoret. Comput. Sci.} 242 (2000) 91--108. 

\bibitem{oJlZ04char} O.~Jenkinson, L.Q.~Zamboni, Characterisations of balanced words via orderings, \emph{Theoret. Comput. Sci} {310} (2004) 247--271. 

\bibitem{jJ00onap} J.~Justin, On a paper by {C}astelli, {M}ignosi, {R}estivo, {\it Theor. Inform. Appl.}  {34} (2000) 373--377.

\bibitem{jJ05epis} J.~Justin, Episturmian morphisms and a {Galois} theorem on continued fractions, 
	{\em Theor. Inform. Appl.} {39} (2005) 207--215.
	
\bibitem{jJgP01frac} J.~Justin, G.~Pirillo, Fractional powers in Sturmian words, \emph{Theoret. Comput. Sci.} {255} (2001) 363--376.	

\bibitem{jJgP02epis} J.~Justin, G.~Pirillo, Episturmian words and episturmian morphisms, \emph{Theoret. Comput. Sci.} {276} (2002) 281--313. 

\bibitem{jJgP02onac} J.~Justin, G.~Pirillo, On a characteristic property 
of Arnoux-Rauzy sequences, {\it Theor. Inform. Appl.} {36} (2002)  
385--388.

\bibitem{jJgP04epis} J.~Justin, G.~Pirillo,  Episturmian words: shifts, morphisms and numeration systems, \emph{Internat. J. Found. Comput. Sci.} {15} (2004) 329--348. 

\bibitem{jJlV00retu} J.~Justin, L.~Vuillon, Return words in {S}turmian and episturmian words, {\em Theor. Inform. Appl.} {34} (2000) 343--356.

\bibitem{tKaV96subs} T.~Komatsu A.J.~van~der~Poorten, Substitution invariant Beatty 
sequences, {\it Jpn. J. Math.} 22 (1996) 349--354.

\bibitem{dK07onst} D.~Krieger, On stabilizers of infinite words, {\it Theoret. Comput. Sci.} 400 (2008) 169--181.

\bibitem{fLgR04quas}
F.~Lev\'{e}, G.~Richomme, Quasiperiodic infinite words: some answers, {\it Bull. Eur. Assoc. Theor. Comput. Sci. (EATCS)} 84 (2004) 128--138.


\bibitem{fLgR07quasB} F.~Lev\'e, G.~Richomme,  Quasiperiodic episturmian words, in: 
 {\em Proceedings of the $6$th International Conference on Words}, Marseille, France, September 17-21, 2007, pp.~201--211.


\bibitem{fLgR07quas} F.~Lev\'{e}, G.~Richomme, Quasiperiodic {S}turmian words and
  morphisms, {\it Theoret. Comput. Sci.} 372 (2007) 15--25.

\bibitem{mL83comb} M.~Lothaire, {\em Combinatorics on Words}, vol.~17 of {\em Encyclopedia of Mathematics and its Applications}, Addison-Wesley, Reading, Massachusetts, 1983.

\bibitem{mL02alge} M.~Lothaire, \emph{Algebraic Combinatorics on
Words}, vol.~90 of {\it Encyclopedia of Mathematics and its Applications}, Cambridge University Press, U.K., 2002.

\bibitem{mL05appl} M.~Lothaire, \emph{Applied Combinatorics on Words}, vol.~105 of {\it Encyclopedia of Mathematics and its Applications}, Cambridge University Press, U.K., 2005.

\bibitem{sM04quas}
S.~Marcus, Quasiperiodic infinite words, {\it Bull. Eur. Assoc. Theor. Comput. Sci. (EATCS)} 82 (2004) 170--174.

\bibitem{fMgP92repe} F.~Mignosi, G.~Pirillo. Repetitions in the Fibonacci 
infinite word, \emph{Theor. Inform. Appl.} {26} (1992) 199--204.

\bibitem{fMpS93morp} F.~Mignosi, P.~S\'{e}\'{e}bold, Morphismes {Sturmiens} et
  r\`{e}gles de {Rauzy}, {\it J. Th\'{e}or. Nombres Bordeaux} 5 (1993) 221--233.



\bibitem{fMlZ02onth} F.~Mignosi, L.Q.~Zamboni, On the number of {A}rnoux-{R}auzy words, 
  {\it Acta Arith.} {101} (2002) 121--129.
  
\bibitem{tM05illu}  T.~Monteil. Illumination dans les billards polygonaux et dynamique symbolique. PhD thesis, Universit\'e de la M\'editerran\'ee, Facult\'e des Sciences de Luminy, December 2005.

\bibitem{tMsM08quas} T.~Monteil, S.~Marcus, Quasiperiodic infinite words: multi-scale case and dynamical properties, {\it Theoret. Comput. Sci.}, to appear, arXiv:math/0603354v1.
  
 \bibitem{gHmM40symb} M.~Morse, G.A.~Hedlund, Symbolic dynamics
II. Sturmian trajectories, \emph{Amer. J. Math.} {62}
(1940) 1--42.
  
\bibitem{gPlV07acha} G.~Paquin, L.~Vuillon, A characterization of balanced episturmian sequences, {\it Electron. J. Combin.} 14 (2007) \#R33, pp.~12. 

\bibitem{gP03ineq} G.~Pirillo, Inequalities characterizing standard 
Sturmian words, {\it Pure Math. Appl.} {14} (2003) 141--144.

\bibitem{gP05ineq} G.~Pirillo, Inequalities characterizing standard 
Sturmian and episturmian words, {\it Theoret. Comput. Sci.} {341} 
(2005) 276--292. 

\bibitem{gP05mors} {G.~Pirillo},  {Morse and {H}edlund's skew {S}turmian words revisited}, {\it Ann. Comb.} 12 (2008) 115--121.

\bibitem{nP02subs} N.~{Pytheas~Fogg}, \emph{Substitutions in Dynamics, Arithmetics and Combinatorics}, vol.~1794 of {\em Lecture Notes in Mathematics}, Springer-Verlag, Berlin, 2002. 

\bibitem{gR82suit} G.~Rauzy, Suites \`a termes dans un alphabet fini, in: {\it S\'emin. Th\'eorie des Nombres}, Exp. No. 25, pp.~16, Univ. Bordeaux I, Talence, 1982--1983.


\bibitem{gR85mots} G.~Rauzy, Mots infinis en arithm\'{e}tique, in: M.~Nivat, D.~Perrin (Eds.), {\it Automata On Infinite Words, Lecture Notes in Computer Science}, vol. 192, Springer-Verlag, Berlin, 1985, pp.~165--171.

\bibitem{gR03conj} G.~Richomme, Conjugacy and episturmian morphisms, \emph{Theoret. Comput. Sci.} {302} (2003) 1--34. 

\bibitem{gR03lynd} G.~Richomme, Lyndon morphisms, {\em Bull. Belg. Math. Soc. Simon Stevin} {10}  (2003) 761--785.


\bibitem{gR07aloc}
G.~Richomme, A local balance property of episturmian
  words, in: {\em Proceedings of the 11th International Conference on Developments in Language Theory 2007} (DLT~'07), July 3--6, Turku, Finland, vol.~4588 of {\em Lecture Notes in Computer
  Science}, Springer, Berlin, 2007, pp.~371--381.


\bibitem{gR07conj}
G.~Richomme, Conjugacy of morphisms and {L}yndon decomposition of
  standard {S}turmian words, {\it Theoret. Comput. Sci.} {380} (2007) 393--400.

\bibitem{gR04onmo} G.~Richomme, On morphisms preserving infinite {L}yndon
  words, {\it Discrete Math. Theor. Comput. Sci.} {9} (2007) 89--108.
  
  \bibitem{gR07priv} G.~Richomme, {\em Private communication}, 2007.

	
\bibitem{rRlZ00agen} R.N.~Risley, L.Q.~Zamboni, 
A generalization of {Sturmian} sequences: {Combinatorial}
structure and transcendence, \emph{Acta Arith.} {95} 
(2000) 167--184.

\bibitem{Si} A.~Siegel, Pure discrete spectrum dynamical systems and periodic tiling associated with a substitution, {\em Ann. Inst. Fourier (Grenoble)} {54} (2004) 341--381.

\bibitem{bTzW06some} B.~Tan, Z.-Y.~Wen, Some properties of the {T}ribonacci sequence, 
{\it European J. Combin.} {28} (2007) 1703--1719.

\bibitem{rT96onco} R.~Tijdeman, On complementary triples of Sturmian bisequences, {\it Indag. Math.} 7 (1996) 419--424.

\bibitem{rT98inte} R.~Tijdeman, Intertwinings of Sturmian sequences, {\it Indag. Math.} 9 (1998) 113--122.

\bibitem{rT00frae} R.~Tijdeman, Fraenkel's conjecture for six sequences, {\it Discrete Math.} 222 (2000) 223--234.

\bibitem{dV00stur} D.~Vandeth, Sturmian words and words with a critical
exponent, \emph{Theoret. Comput. Sci.} {242} (2000) 283--300.

\bibitem{pV86symb} P.~Veerman, Symbolic dynamics and rotation numbers, {\em Physica A} {134} (1986) 543--576.

\bibitem{pV87symb} P.~Veerman, Symbolic dynamics of order-preserving orbits, {\em Physica D} {29} (1987) 191--201.


\bibitem{lV01acha} L.~Vuillon, A characterization of {S}turmian words by return words, 
  {\it European J. Combin.} 22 (2001) 263--275.

\bibitem{lv03bala} L.~Vuillon, Balanced words, {\it Bull. Belg. Math. Soc. Simon Stevin} 10 (2003) 787--805.

\bibitem{zWyZ99some} Z.-X.~Wen, Y.~Zhang, Some remarks on invertible substitutions on three letter alphabet, {\em Chinese Sci. Bull.} 44 (19) (1999) 1755--1760. 

\bibitem{nWlZ01freq} N.~Wozny, L.Q.~Zamboni,  Frequencies of factors in
  {A}rnoux-{R}auzy sequences, {\it Acta Arith.} {96} (2001) 261--278.
  
\bibitem{sY99onst}  S.-I.~Yasutomi, On Sturmian sequences which are invariant under some 
substitutions, in: {\em Number theory and its applications} (Kyoto, 1997), Kluwer Acad. Publ., Dordrecht, 1999,  pp.~347--373.
  
\bibitem{lZ98uneg} L.Q.~Zamboni, Une g\'en\'eralisation du th\'eor\`eme de Lagrange sur le d\'eveloppement en fraction continue, {\it C.R. Acad. Sci. Paris S\'er. I Math.} {327} (1998) 527--530. 


\end{thebibliography}
\end{document}